\documentclass{article}

\usepackage{amsmath}
\usepackage{amsfonts}

\newtheorem{theo}{Theorem}

\newtheorem{prop}{Proposition}
\newtheorem{lemm}{Lemma}
\newtheorem{defn}{Definition}

\def\N{\mathbb{N}}
\def\E{\mathbb{E}}
\def\0{{\bf 0}}

\def\R{\mathbb{R}}

\def\dist{ {\rm dist}}
\def\B{{\cal B}}

\def\S{{\cal S}}

\renewcommand{\E}{\mathbb E \,}
\newcommand{\C}{{\cal C}}

\newcommand{\tod}{\stackrel{{\cal D}}{\longrightarrow}}
\newcommand{\eqd}{\stackrel{{\cal D}}{=}}

\newcommand{\eqco}{\setcounter{equation}{0}}
\newcommand{\thco}{\setcounter{theo}{0}}
\newcommand{\prco}{\setcounter{prop}{0}}
\newcommand{\laco}{\setcounter{lemm}{0}}
\newcommand{\coco}{\setcounter{coro}{0}}
\newcommand{\cjco}{\setcounter{conj}{0}}

\newcommand{\deco}{\setcounter{defn}{0}}
\newcommand{\allco}{\eqco  \thco \prco \laco \coco \cjco \deco}
\newcommand{\qed}{\hfill{\rule[-.2mm]{3mm}{3mm}}}

\newcommand{\Po}{{\cal P}}
\newcommand{\X}{{\cal X}}

\renewcommand{\H}{{\cal H}}
\renewcommand{\P}{{{\cal P}}}

\newcommand{\A}{{\cal A}}
\newcommand{\Cov}{{\rm Cov}}

\newcommand{\Var}{{\rm Var}}
\newcommand{\var}{{\rm Var}}
\newcommand{\An}{{\rm An}}

\newcommand{\diam}{{\rm diam}}

\newcommand{\K}{{\cal K}}

\newcommand{\sr}{{\text{Poisson-like}}}

\newcommand{\tH}{\tilde{\cal H}}

\def\bdm{\begin{displaymath}}
\newcommand{\edm}{\end{displaymath}}
\def\benu{\begin{enumerate}}
\def\eenu{\end{enumerate}}
\def\beqn{\begin{equation}}
\def\eeqn{\end{equation}}
\def\be{\begin{equation}}
\def\ee{\end{equation}}
\def\bea{\begin{eqnarray}}
\def\eea{\end{eqnarray}}
\newcommand{\bean}{\begin{eqnarray*}}
\newcommand{\eean}{\end{eqnarray*}}
\newcommand{\bear}{\begin{eqnarray}}
\newcommand{\eear}{\end{eqnarray}}

\renewcommand{\epsilon}{\varepsilon}

\def\Comment#1{
\marginpar{$\bullet$\quad{\tiny #1}}}

\def\R{\mathbb{R}}

\def\vol{\hbox{\rm vol}}

\def\de{{\delta}}
\def\A{{\cal A}}

\def\Comment#1{\lineskip-4pt
\marginpar{ $\bullet$\quad{\em\small #1}}}

\def\qed{\hfill\hbox{${\vcenter{\vbox{
    \hrule height 0.4pt\hbox{\vrule width 0.4pt height 6pt
    \kern5pt\vrule width 0.4pt}\hrule height 0.4pt}}}$}}

\def\la{{\lambda}}
\def\ka{{\kappa}}

\begin{document}

\title{\bf Stabilization and limit theorems for geometric functionals
           of Gibbs point processes}

\author{T. Schreiber and J. E. Yukich$^{1}$ }

\date{\today}
\maketitle

\footnotetext{ {\em American Mathematical Society 2000 subject
classifications.} Primary 60F05, 60G55, Secondary 60D05}

 \footnotetext{
{\em Key words and phrases.} Gibbs point processes, Gaussian
limits, hard core model, random packing, geometric graphs, loss
networks}


\footnotetext{$~^1$ This material is based upon work supported by
the NSA under  grant H98230-06-1-0052}

\begin{abstract}
 Given a Gibbs point
process $\P^{\Psi}$ on $\R^d$ having  a weak enough potential $\Psi$, we
consider the random measures $\mu_\la := \sum_{x \in \P^{\Psi}
\cap Q_\la} \xi(x, \P^{\Psi} \cap Q_\la) \delta_{x/\la^{1/d}}$,
where $Q_{\la} := [-\la^{1/d}/2,\la^{1/d}/2]^d$ is the volume
$\la$ cube and where $\xi(\cdot,\cdot)$ is a translation invariant
stabilizing functional.
Subject to $\Psi$ satisfying a localization property and
translation invariance, we establish weak laws of large numbers
 for
$\la^{-1} \mu_\la(f)$, $f$ a bounded test function on $\R^d$, and
weak convergence of $\la^{-1/2} \mu_\la(f),$ suitably centered, to
a Gaussian field acting on bounded test functions. The result
yields limit laws for geometric functionals on Gibbs point
processes including the Strauss and area interaction point
processes as well as more general point processes defined by the
Widom-Rowlinson and hard-core model. We provide applications to
random sequential packing on Gibbsian input, to functionals of
Euclidean graphs, networks, and percolation models on Gibbsian
input, and to quantization via Gibbsian input. 

\end{abstract}

\section{Introduction}\label{Introduction}

\allco


 Functionals of large complex geometric structures  often consist
of sums of spatially dependent terms admitting the representation
\be \label{sums} \sum_{x \in \X} \xi(x,\X), \ee where $\X \subset
\R^d$ is locally finite and where the function $\xi$, defined on
all pairs $(x, \X)$, with $x \in \X$, represents the interaction
of $x$ with respect to $\X$. When $\X$ is a random $n$ point set
in $\R^d$ (i.e. a finite spatial point process), the asymptotic
analysis of the suitably scaled sums (\ref{sums})  as $n \to
\infty$
 can often be handled by  $M$-dependent
methods, ergodic theory, or mixing methods. However there are
situations where these classical methods are either not directly
applicable,  do not give explicit asymptotics in terms of
underlying geometry and point densities, or do not easily yield
explicit rates of convergence.   Stabilization methods originating
in \cite{PY1} and further developed in \cite{BY2, PY2, PY5},
provide another approach for handling sums of spatially dependent
terms.

There are several similar definitions of stabilization, but the
essence is captured by the notion of stabilization of the
functional $\xi$ with respect to a rate $\tau > 0$ homogeneous Poisson
point process $\P := \P_\tau$ on $\R^d$, defined as follows. Say that
$\xi$ is translation invariant if $\xi(x,\X) = \xi(x + z, \X + z)$
for all $z \in \R^d$. Let $B_r(x)$ denote the Euclidean ball
centered at $x$ with radius $r \in \R_+ := [0, \infty)$. Letting
$\0$ denote the origin of $\R^d$, we say that a translation
invariant $\xi$ is {\em stabilizing} on $\P = \P_\tau$
if there exists
an a.s. finite random variable $R:=R^{\xi}(\tau)$ (a `radius of
stabilization') such that
 \be \label{stab}
\xi(\0, \P \cap B_R(\0)) = \xi(\0, \P \cap B_R(x)
\cup \A ) \ee
 for all locally finite $\A \subset \R^d \setminus B_R(\0)$.



 Consider the point measures \be \label{rm} \mu_{\la }
:= \sum_{x \in \P \cap Q_\la} \xi(x, \P \cap
Q_\la)\delta_{x/\la^{1/d}}, \ee where $\delta_x$ denotes the unit
Dirac point mass at $x$ whereas $Q_{\la} := [-\la^{1/d}/2,\la^{1/d}/2]^d$
is the $\la$-volume cube. Let $\B(Q_1)$ denote the class of all
bounded $f: Q_1 \to \R$ 
and for all random
point measures $\mu$ on $\R^d$ let $\langle f, \mu \rangle:= \int
f d\mu$ and let $\bar{\mu}: = \mu - \E[\mu]$.

Stabilization of translation invariant $\xi$ on $\P$, as
defined in (\ref{stab}), together with  stabilization of $\xi$ on
$\P \cap Q_\la, \la \geq 1$, 
 when
combined with appropriate moment conditions on $\xi$, yields for
all $f \in \B(Q_1)$ the law of large numbers \cite{Pe2, PY4}
 \be
\label{WLLNa} 
\lim_{\la \to \infty}  \la^{-1}\langle f, \mu_{\la} \rangle = \tau
\E [\xi(\0, \P)]  \int_{Q_1} f(x)  dx \ \ \text{in} \ \
L^1 \ \ \text{and in} \ \ L^2,\ee and, if the stabilization radii
on $\P$ and $\P \cap Q_\la, \la \geq 1,$ 
  decay exponentially fast,
  then  \cite{BY2, Pe1}
 \be
\label{VARa} \lim_{\la \to \infty} \la^{-1} \Var [\langle f,
\mu_{\la} \rangle] = \tau V^{\xi}(\tau ) \int_{Q_1} f(x)
 dx,
 \ee
where  for all  $\tau> 0$
$$
V^{\xi}(\tau):= \E[ \xi(\0, \P)^2] + \tau \int_{\R^d} [\E
\xi(\0,\P \cup \{z\}) \xi(z, \P \cup \{\0\}) - \E[
\xi(\0, \P)^2]]dz.
$$
Additionally, the   finite-dimensional distributions $(\langle
f_1, \la^{-1/2} \overline{\mu}_{\la}^{\xi} \rangle, \ldots,
 \langle f_k, \la^{-1/2}\overline{\mu}_{\la}^{\xi} \rangle),$
$f_1,\ldots,f_k \in \B(Q_1),$ converge to a Gaussian field with
covariance kernel
 \be
\label{covkernela} (f,g) \mapsto  \tau V^{\xi}(\tau ) \int_{Q_1}
f(x)g(x)  dx.
  \ee

The limits (\ref{WLLNa})-(\ref{covkernela}) establish asymptotics
for functionals and measures defined in terms of independent input
and one might expect analogous asymptotics  for functionals  of
dependent input subject to  weak long range dependence conditions.
The main purpose of this paper is to show that this is indeed the
case.  We establish the analogs of
(\ref{WLLNa})-(\ref{covkernela}) when $\P = \P_\tau$ is replaced
by a weak Gibbsian
 modification having an exponentially localized potential; see
 Theorems \ref{WLLN}-\ref{CLT} for  a precise statement of the limit theory
 for functionals of Gibbsian input.
Gibbsian point processes covered by this generalization include,
for low enough reference intensity $\tau$,  the Strauss process,
the area interaction process, as well as point processes defined
by the continuum Widom Rowlinson and hard-core models.
 Gibbsian point
processes considered here are intrinsically algorithmic. Their
computational efficiency yields numerical estimates for asymptotic
limits appearing in our main results.

Functionals of geometric graphs over Gibbsian input on large
cubes, as well as functionals of  random sequential packing models
defined by Gibbsian input on large cubes, consequently satisfy
weak laws of large numbers and central limit theorems as the cube
size tends to infinity. The precise limit theorems are provided in
sections 6 and 7, which also includes asymptotics for functionals
of communication networks and continuum percolation models over
Gibbsian point sets, as well as asymptotics for the distortion
error arising in Gibbsian quantization of probability measures.

\vskip.5cm


\section{Gibbs point processes and their stabilizing functionals}

\allco

\subsection{Gibbs point processes with localized potentials}\label{KONSGRAPH}

 Throughout $\Psi$ denotes a translation invariant functional  defined on locally finite
 collections of points $\X$ in ${\Bbb R}^d$
 and admitting values in ${\Bbb R}_+ \cup \{ + \infty \}.$
By translation invariant we mean $\Psi(\X) = \Psi(y + \X)$ for all
$y \in \R^d$.
 In the sequel we refer to $\Psi$
  as the potential, Hamiltonian or energy functional, with all three terms
 used interchangeably.
 For a locally finite 
  point set ${\cal X}$ in
 ${\Bbb R}^d$ and
 an open bounded set $D \subseteq {\Bbb R}^d$ we define $\Psi_D({\cal X}) :=
 \Psi({\cal X} \cap D).$   We shall always
 assume that for all open, bounded $D$ the potential $\Psi_D({\cal P})$ admits finite values
 with non-zero probability, where we recall that ${\cal P}:= {\cal P}_{\tau}$ is the
 Poisson point process of some arbitrary but fixed intensity $\tau > 0$
 in ${\Bbb R}^d$.  Moreover, we assume the
 Hamiltonian is {\it hereditary}, 
 that is to say if $\Psi({\cal X}) = +\infty$
 for some ${\cal X}$ then $\Psi({\cal Y}) = +\infty$ for all
 ${\cal Y} \supseteq {\cal X}.$ This puts us
 in a position to define the Gibbs point process ${\cal P}^{\Psi}_D$
 given in law by
 \begin{equation}\label{GIBBS1}
  \frac{d{\cal L}({\cal P}^{\Psi}_D)}{d{\cal L}({\cal P})}({\cal X}) :=
  \frac{\exp (-\Psi_D({\cal X}) )}{Z[\Psi_D]},
 \end{equation}
 where $ Z[\Psi_D] := {\Bbb E}\exp(-\Psi_D({\cal P}))$
 is the normalizing constant for (\ref{GIBBS1}), also called the
 partition function.

The following definition is central to this paper.

\begin{defn} \label{centraldef}
 For a decreasing right-continuous
 function $\psi : {\Bbb R}_+ \to [0,1]$ with $\lim_{r\to\infty} \psi(r)=0$
 we say that the  Hamiltonian $\Psi$ is $\psi$-localized
 iff for each $x \in {\Bbb R}^d,$ each finite ${\cal X}$ and each $r > 0$ the
 {\it add-one potential}, inheriting from $\Psi$ its translation invariance
 property,
 $$ \Delta(x,{\cal X}) := \Psi({\cal X} \cup \{ x \}) - \Psi({\cal X}), $$
 satisfies
  \begin{equation}\label{LOCALIZATION}
  0 \leq \Delta_{[r]}(x,{\cal X} \cap B_r(x)) \leq \Delta(x,{\cal X}) \leq
  \Delta^{[r]}(x,{\cal X} \cap B_r(x)),
 \end{equation}
 where $\Delta_{[r]}(\cdot,\cdot)$ and $\Delta^{[r]}(\cdot,\cdot)$ are certain
 translation invariant deterministic functionals such that
  uniformly in $x \in \R^d$ and $\X \subset \R^d$ we have
 \begin{equation}\label{LCBOUND}
  0 \leq \exp(-\Delta_{[r]}(x,{\cal X} \cap B_r(x))) -
         \exp(-\Delta^{[r]}(x,{\cal X} \cap B_r(x))) \leq \psi(r).
 \end{equation}
\end{defn}

 In other words, even though determining exactly the value of the add-one
 potential $\Delta(x,{\cal X})$ may require the knowledge of the whole
 configuration ${\cal X},$ knowing just ${\cal X} \cap B_r(x)$ we can
 determine the value of $\exp(-\Delta(x,{\cal X}))$ with accuracy at least
 $\psi(r)$ which tends to $0$ as $r\to\infty.$ In case where both
 $\Psi({\cal X} \cup \{ x \})$ and $\Psi({\cal X})$ are $+\infty$ we
 set by convention $\Delta(x,{\cal X}) := 0.$ We also require that
 $ \Delta(x,\emptyset) < + \infty$
 to prevent the Gibbs process  ${\cal P}^{\Psi}_D$ from concentrating
 on $\emptyset.$
 The functionals $\Delta_{[r]}(\cdot,\cdot)$ and $\Delta^{[r]}(\cdot,\cdot)$
 will be called {\em lower and upper add-one potentials} respectively.
 Note that
 the required non-negativity of the add-one potential is not  particularly
 restrictive  because whenever the add-one potential admits a finite
 lower bound, possibly negative $-a < 0,$ it can be reduced to the present
 setting by adding $a |{\cal X}|$ to $\Psi$ and by replacing the underlying intensity
 $\tau$   with $\tau \exp(a).$ Imposing the presence of a lower
 bound for the add-one potential or other related growth conditions is a
 usual assumption to avoid density explosions and infinite values of the
 partition function in (\ref{GIBBS1}), see \cite{Rue70}. 

 Every Poisson point process has a
$\psi$-localized potential, since in this case $\Psi \equiv 0$ and
thus $\Delta \equiv 0$. Less trivially, a large number of Gibbs
point processes, including those in modelling problems in
statistical mechanics, communication networks, and biology have
$\psi$-localized potentials. This includes the Strauss process,
the area interaction process, processes having finite and infinite
range pair potential functions, and the hard-core and
Widom-Rowlinson models; see section \ref{S41} for details.

\subsection{Graphical construction of Gibbs point processes with localized potentials}

\label{graphconstruct}


 For a $\psi$-localized potential $\Psi$ the resulting Gibbs point process
 ${\cal P}^{\Psi}_D$ admits a particularly convenient graphical construction
 in the spirit of Fern\'andez, Ferrari and Garcia \cite{FFG1}-\cite{FFG3}.
 While adding a number of new ideas, in our presentation below we follow
\cite{FFG1}-\cite{FFG3} as well as the developments
 in \cite{BESY}.
 Consider a stationary homogeneous {\it free birth and death process}
 $(\rho^D_t)_{t \in {\Bbb R}}$ in $D$ with the following dynamics:
 \begin{itemize}
  \item A new point $x \in D$ is born in $\rho^D_t$ during the time interval $[t-dt,t]$
        with probability $\tau dx dt,$
  \item An existing point $x \in \rho^D_t$ dies during the time interval
        $[t-dt,t]$ with probability $dt,$ that is
        the lifetimes of points of the process are independent standard
        exponential.
 \end{itemize}
 Clearly, the unique stationary and reversible measure for this process is just
 the law of the Poisson point process ${\cal P} \cap D.$

 Consider now the
 following {\it trimming} procedure performed on $\rho^D_t,$ based on
 the ideas developed in \cite{FFG1}-\cite{FFG3}. Choose a birth
 site for a point $x \in D$ at some time $t \in \R$ and draw a random number
 $\eta \in {\Bbb R}_+$ from the law given by the distribution function
 $1-\psi(\cdot).$ Then, accept it with probability $\exp(-\Delta^{[r]}(x,\rho^D_{t-}
 \cap B_{\eta}(x)))/\psi(\eta)$ and reject with the complementary probability
 if the acceptance/rejection statuses of all points in $\rho^D_{t-} \cap B_{\eta}(x)$
 are determined, otherwise proceed recursively to determine the statuses of
 points in $B_{\eta}(x).$

 Before discussing any further properties of this
 procedure, we have to ensure first that it actually  terminates. To this end,
 note that each point $x$ with the property of having the ball $B_{\eta}(x)$
 devoid of points from $\rho^D_{t-}$ at its birth time $t$ has its acceptance
 status determined. More generally, the acceptance status of a point $x$
 at its birth time $t$ only depends on the status of points in
 $\rho^D_{t-} \cap B_{\eta}(x),$ that is to say points born before and
 falling into $B_{\eta}(x).$ We call these points {\it causal ancestors}
 of $x$ and, in general, for a subset $A \subseteq D$ by $\An_t[A]$
 we denote the set of all points in $\rho^D_t \cap A,$ their causal
 ancestors, the causal ancestors of their ancestors and so forth
 throughout all past generations. The set $\An_t[A]$ is referred to
 as the {\it causal ancestor cone} or {\it causal ancestor clan} of $A$
 with respect to the birth and death process $(\rho^D_t)_{t \in {\Bbb
 R}}$.

 It is now clear that in order for our recursive status determination
 procedure to terminate for all points of $\rho^D_t$ in $A,$ it is
 enough to have the causal ancestor cone $\An_t[A]$ finite. This
 is easily checked to be a.s. the case for each $A \subseteq D$ \---
 indeed, since $D$ is bounded, a.s. there exists some $s < t$ such
 that $\rho^D_s = \emptyset$ and thus no ancestor clan of a point
 alive at time $t$ can go past $s$ backwards in time.

 Having defined the trimming procedure above, we recursively remove
 from $\rho^D_t$ the points rejected at their birth, and we write
 $(\gamma^D_t)_{t\in {\Bbb R}}$ for the resulting process. Clearly,
 $\gamma^D_t$ is stationary because so was $\rho^D_t$ and the
 acceptance/rejection procedure is time-invariant as well.
 Moreover, the process $\gamma^D_t$ is easily seen to evolve
 according to the following dynamics:
 \begin{itemize}
  \item Add a new point $x$  with intensity $\tau \exp(-\Delta(x,\gamma^D_t)) dx dt,$
  \item Remove an existing point with intensity $dt.$
 \end{itemize}
 These are the standard Monte-Carlo dynamics for ${\cal P}^{\Psi}_D \cap D$
 as given in (\ref{GIBBS1}) and the law of ${\cal P}^{\Psi}_D$ is its unique
 invariant distribution. Consequently, in full analogy with \cite{FFG1}-\cite{FFG3}
 the point process $\gamma^D_t$ coincides in law with ${\cal P}^{\Psi}_D$ for
 all $t \in {\Bbb R}.$

\subsection{Exponentially localized potentials and infinite volume limits}

Recalling  the definition of $\psi$ from Definition
\ref{centraldef}, we henceforth  assume that there is a $C_1
> 0$ such that
 \begin{equation}\label{ExpDec}
  \psi(r) \leq \exp(-C_1r) \ \ \forall r > 0.
 \end{equation}
 It should be emphasized  that we require (\ref{ExpDec}) to hold
 {\it for all} $r>0$ and not just for $r$ large enough.
 It is known, see \cite{FFG1}-\cite{FFG3} where a proof
 based on subcritical branching process domination is given, that if $C_1$
 is chosen large enough, then all causal ancestor cones are a.s. finite and,
 in fact, there is a $C_2 > 0$ such that for all $t, R \in \R_+$
 and $A \subset D$ we have the crucial bound
 \begin{equation}\label{ExpDec2}
  {\Bbb P} [\diam \An_t[A] \geq R + \diam(A)] \leq \vol(A) \exp(-C_2
  R).
 \end{equation}
 Moreover, the constant $C_2$ in (\ref{ExpDec2})  does not depend on $D.$
 If (\ref{ExpDec}) is satisfied with the constant $C_1$ large enough so
 that
 (\ref{ExpDec2}) holds as well, then the potential $\Psi$ is declared
 {\it exponentially localized}. 


 Putting $D_n:= [-n,n]^d$, this puts us in a
 position to construct the infinite volume limit (thermodynamic limit) for ${\cal
 P}^{\Psi}_{D_n}$ as $n \to \infty$.
 Indeed, consider the
 infinite volume version $\rho_t$ of our stationary free birth and death
 process $\rho^D_t,$ constructed as $\rho^D_t$ with $D$ replaced by
 $\R^d$.
 Clearly, for each $t \in {\Bbb R}$ we have that $\rho_t$ coincides in law
 with ${\cal P}.$
 Moreover, in view of (\ref{ExpDec2}) and recalling
 that $c$ there did not depend on $D,$ we see that the trimming procedure
 as described above is also valid for the infinite volume process $\rho_t,$
 yielding the stationary {\it trimmed} process $\gamma_t.$
These remarks justify defining the following point process, used
in all that follows.

\begin{defn} We define the
 thermodynamic limit ${\cal P}^{\Psi} := \P_\tau^{\Psi}$ to be the point process coinciding
 in law with $\gamma_0$ and hence with $\gamma_t$ for all $t.$
\end{defn}

 To provide some further motivation for granting to ${\cal P}^{\Psi}$
 the name of {\it thermodynamic limit} note that
 the process ${\cal P}^{\Psi}$ enjoys the following important property:
 for any bounded set $D \subseteq {\Bbb R}^d$, any locally finite
 point configuration ${\cal X} \subseteq D^c$ and any finite point configuration
 ${\cal Y} \subseteq D$ the conditional law of
 ${\cal P}^{\Psi} \cap D$ on the event ${\cal P}^{\Psi} \cap D^c =
  {\cal X}$ is given by
 \begin{equation}\label{CondSpec}
    \frac{d{\cal L}({\cal P}^{\Psi} \cap D | {\cal P}^{\Psi} \cap D^c = {\cal X})}
         {d{\cal L}({\cal P} \cap D)}[{\cal Y}] =
    \frac{\exp(-\Psi({\cal Y}|{\cal X}))}{Z_D[\Psi|{\cal X}]},
 \end{equation}
 where
 $$ Z_D[\Psi|{\cal X}] = {\Bbb E} \exp(-\Psi({\cal P} \cap D |{\cal X})) $$
 whereas
 $$ \Psi({\cal Y}|{\cal X}) := \lim_{r \to \infty} \Psi({\cal Y} \cup {\cal X} \cap B(\0,r))
    - \Psi({\cal X} \cap B(\0,r)) $$
 with the existence of the limit guaranteed by the localization
 condition (\ref{LOCALIZATION}). Moreover, the so-constructed ${\cal P}^{\Psi}$ is the only point process
 with the above properties \--- to see it take $D_n := [-n,n]^d \uparrow {\Bbb R}^d$
 and note that in view of the graphical construction specialized for the conditional
 specification (\ref{CondSpec}), the relation (\ref{ExpDec2}) guarantees that the
 process in some fixed bounded $A \subseteq \R^d$ exhibits exponentially
 decaying dependencies on the external configuration in $D_n^c$ as $n\to\infty.$

\subsection{Stabilizing functionals of Gibbs point processes}

 In this section we specialize to our Gibbs point process setting
 the notion of a stabilizing functional, see \cite{BY2,PY1,PY2,PY4} and the
 references therein.  As in section 1, let  $\xi(\cdot,\cdot)$ be a translation invariant functional
 defined on pairs $(x,{\cal X})$ where ${\cal X}$ is a finite point collection
 in ${\Bbb R}^d$ and $x \in {\cal X}$.
 Further, when $x \notin \X$, we abbreviate $\xi(x, \X \cup \{x\})$ by $\xi(x,\X)$.

Next, suppose that a given point process $\Xi$ on ${\Bbb R}^d$ is
stochastically dominated by a homogeneous Poisson point process
and suppose that there exists $C_3 > 0$ such that for every ball
$B_r(x)$ the conditional probability of $B_r(x)$ not being hit by
$\Xi$ given the external configuration
${\cal E}:= \Xi \setminus B_r(x)$ admits the bound \be \label{LBd}
    {\Bbb P}[\Xi \cap B_r(x) = \emptyset|{\cal E}] \leq \exp(-C_3 r^d)
 \ee
uniformly in ${\cal E}.$
Stochastic domination and (\ref{LBd}) provide upper and lower
stochastic bounds on the number of points in any ball analogous to
those satisfied by a homogeneous Poisson point process and for
this reason such $\Xi$ are called  {\em ${\sr}$}. The next
proposition tells us that the Gibbs point processes considered
here are $\sr$.

\begin{prop} Every Gibbs point process $\P^{\Psi}$ with an
exponentially localized potential is a $\sr$ process.
\end{prop}

{\em Proof.} Indeed, the stochastic domination by ${\cal P}$ comes
from the obvious relation $\gamma_0 \subseteq \rho_0$ in the above
graphical construction of ${\cal P}^{\Psi}$ because $\rho_0$
coincides in law with ${\cal P}.$ The second relation (\ref{LBd})
follows by the graphical construction as well. Indeed, we have
$\Delta(x,\emptyset) < \infty$ and hence, by (\ref{LOCALIZATION})
and (\ref{LCBOUND}), in the course of the dynamics given by the
graphical construction the acceptance probability for a birth
attempt at some $y$ inside a ball $B_{r-s}(x)$ with no points
alive in the whole $B_r(x)$ is uniformly bounded away from $0,$
both in the location of the point $y$ attempting to be born and in
the external configuration,
as soon as $s$ and $r>s$
are taken large enough. On the other hand, the ball reaches a
completely empty state with intensity at most  $1.$ Consequently,
the time fraction of having the ball fully empty decays
exponentially with the volume of the ball uniformly in the
external configuration and hence so does the probability of having
no point alive in $B_r(x)$ at the time $0$ by stationarity of the
graphical construction in time. \qed


\ \

 Similarly to (\ref{stab}),  say that $\xi$ is a
 {\it stabilizing functional in the wide sense} if for every $\sr$ process
 $\Xi$ there exists an a.s. finite stabilization radius $R := R^{\xi}(x,\Xi),$  such that a.s.
 \begin{equation}\label{StabRelW}
  \xi(x,\Xi \cap B_R(x)) = \xi(x,[\Xi \cap B_R(x)] \cup {\cal A})
 \end{equation}
 for all locally finite point collections ${\cal A} \subseteq \R^d \setminus B_R(x).$
 Stabilizing functionals in the wide sense can  a.s. be extended to the whole process
 $\Xi$, that is to say for all $x \in \R^d$
 $$ \xi(x,\Xi) := \lim_{r\to\infty} \xi(x,\Xi\cap B_R(x)) $$
 is a.s. well defined.

Given $s > 0$ and a Poisson-like process $\Xi$ define the tail
probability
$$
\tau(s):= \tau(s, \Xi):= \max \left[ \sup_{\la \geq 1, x \in
Q_\la} {\Bbb P} [R(x, \Xi \cap Q_\la) > s], \ \ {\Bbb P} [R(x, \Xi
) > s] \right].
$$

 Further, we say that $\xi$ is {\it exponentially stabilizing in the wide sense}
 if for every $\sr$ process $\Xi$
 we have $\limsup_{s \to \infty} s^{-1} \log \tau(s) < 0$. Thus, if
 $\xi$ is exponentially stabilizing in the wide sense, then
 there exists a $C_4$ such that for all $s \in \R_+$ we have
   \begin{equation}\label{ExpStab}
   \sup_{x \in \R^d}
    {\Bbb P}[R(x,\Xi) > s] \leq \exp(-C_4 s) \
    \text{and} \ \
    \sup_{\la \geq 1, x \in Q_\la} {\Bbb P} [R(x, \Xi \cap Q_\la) > s]
    \leq \exp(-C_4 s).
 \end{equation}

We stress that, unlike in the standard Poisson input setting of
\cite{BY2,PY1,PY2,PY4}, where Poisson points in disjoint sets are
independent, the configuration $\Xi \cap B_R(x)$ will, in general,
depend on the configuration in $\Xi \cap B_R(x)^c$. Thus, unlike
the standard Poisson input setting, the wide sense stabilization
of $\xi$ at $x$ within radius $R$ does not imply that the value of
$\xi(x,\Xi)$ does not depend on the configuration outside
$B_R(x);$ on the other hand this value is independent of the
configuration in $B_R(x)^c$ {\it given the configuration} $\Xi
\cap B_R(x)$. This weak dependence feature of wide sense
stabilization, which carries additional technical considerations,
allows  us to establish limit theory for functionals and measures
in geometric probability over point sets more general than the
usual Poisson and binomial point sets.

 As we will see shortly,  many functionals which stabilize in the
 standard Poisson input setting also stabilize in the wide sense.
 Possibly there are some functionals  which
 stabilize over Poisson samples but which do stabilize in the wide sense, but we
 are not aware of these functionals.
 For these reasons,  when the context is clear, {\em we  will henceforth abuse terminology and use the
 term `stabilization' to mean  `stabilization
 in the wide sense'}, with a similar meaning for `exponentially stabilizing'.

\subsection{Functionals with bounded perturbations}
 The theory presented in this paper is mainly confined to {\it translation invariant}
 geometric functionals and its extension to non-translation invariant
 functionals seems to require non-trivial effort. Nevertheless, a small step
 towards the non-translation invariant set-up can be made with only slight
 modification of the existing theory. This extension is the subject of the
 present subsection and it deals with asymptotically negligible bounded
 perturbations  of translation-invariant functionals. To put it in formal terms,
 consider the following notion. Consider a \sr \ input point process $\Xi.$
 Assume that $\xi(\cdot,\cdot)$ is a
 translation invariant geometric functional exponentially stabilizing in the
 wide sense and let $\hat\xi(\cdot,\cdot;\la)$ be a family of geometric
 functionals indexed by the extra parameter $\la > 0,$ not assumed to be
 translation invariant but enjoying the following properties:
 \begin{itemize}
  \item For each $\la > 0$ the functional $\hat\xi(\cdot,\cdot;\la)$ admits
        a representation
        \begin{equation}\label{BDDIST1}
         \hat\xi(x,\X;\la) = \xi(x,\X) + \delta(x,\X;\lambda),
        \end{equation}
        where the correction (perturbation) $\delta(x,\X;\la)$ is not
        necessarily translation invariant but, for all $p > 0$
        it satisfies  the moment bound
        \begin{equation}\label{BDDIST2}
         \sup_{x \in \R^d} \E [ \delta(x,\Xi;\lambda)]^p \leq
         \epsilon(\la,p) < \infty,
        \end{equation}
        where $\lim_{\la \to \infty} \epsilon(\la,p) = 0$ for each fixed $p.$
  \item The perturbation  $\delta(\cdot,\cdot;\la)$ satisfies the
        wide sense exponential stabilization with the same stabilization
        radius $R^{\xi}(\cdot,\cdot)$ as $\xi.$
  \end{itemize}
  If these two conditions hold, we say that $\hat\xi(\cdot,\cdot;\la)$ is
  an {\em asymptotically negligible bounded perturbation of $\xi$ on input
  $\Xi,$}; for brevity we call it just a bounded perturbation  of $\xi$
  in the sequel. The message of
  this subsection, to be made formal below, is that 
  the asymptotic behavior of bounded perturbations  of a
  translation invariant functional
  is indistinguishable from the asymptotic properties of the functional
  itself. This observation brings the limit theory for stochastic quantization within the compass of
  stabilizing functionals; see section \ref{QUANT}.

\section{Weak laws of large numbers and central limit theorems}

\allco We now state our main results, which show that sums of
 stabilizing functionals defined on Gibbsian input
 (with exponentially localized potential) on large cubes satisfy
 weak laws of large numbers and Gaussian limits
as the cube size tends to infinity.  For all $\la > 0$, let
$Q_{\la} := [-\la^{1/d}/2,\la^{1/d}/2]^d$ be the volume $\la$ cube
 centered at the origin of $\R^d$,
and let $\mu^{\xi}_{\la}$ be the $\la$-rescaled $\xi$-empirical
measure on
 $Q_1 := [-1/2,1/2]^d$, that is
 \be \label{measure} \mu^{\xi}_{\la} := \sum_{x \in {\cal P}^{\Psi} \cap Q_{\la}} \xi(x,{\cal P}^{\Psi} \cap Q_{\la})
    \delta_{x/\la^{1/d}}. \ee
Let  $H^{\xi}_{\la}:= \mu_{\la}^{\xi}(Q_1)$ be the total mass of
$\mu^{\xi}_{\la},$ and for future reference,  define also the
non-rescaled infinite-volume measure
 \begin{equation}\label{MUX}
  \mu^{\xi} := \sum_{x \in {\cal P}^{\Psi}} \xi(x,{\cal P}^{\Psi}) \delta_x.
 \end{equation}
 Let $p \in [0, \infty)$.
Say that $\xi$ satisfies the $p$-moment condition if \be
\label{mom}\sup_\la \sup_{x \in {\cal P}^{\Psi} \cap Q_{\la}, \ \X
\in \C } \E [|\xi(x,{\cal P}^{\Psi}  \cup \X)|^p] < \infty, \ee
where $\C$ denotes the collection of all finite point sets in
$\R^d$.

Recall that $\B(Q_1)$ denotes the set of bounded
  $f:Q_1 \to {\Bbb R}$ and that $\bar{\mu}^{\xi}_{\la}:= \mu^{\xi}_{\la} -
\E[\mu^{\xi}_{\la}]$.  Under appropriate moment conditions, our
first two results establish a weak law of large numbers and
variance asymptotics for $\langle f, \mu^{\xi}_{\la} \rangle, \ f
\in \B(Q_1)$, as $\la \to \infty$. Our third result shows that the
finite-dimensional distributions of $( \la^{-1/2}\langle f_1,
\bar{\mu}^{\xi}_{\la} \rangle,...,\la^{-1/2}\langle f_m,
\bar{\mu}^{\xi}_{\la} \rangle), \ f_1,...,f_m \in \B(Q_1),$
converge to those of a multivariate normal as $\la \to \infty$,
and, in the univariate CLT we establish a rate of convergence.
Finally our last general result establishes asymptotics for
bounded perturbations  of a translation invariant $\xi$.


\begin{theo}{(WLLN)}\label{WLLN}
Assume that $\xi$ is stabilizing and satisfies the $p$-moment
condition (\ref{mom}) for some $p > 1$.
We have for each $f \in \B(Q_1)$ \be \label{wllnlimit}
  \lim_{\la\to\infty} \la^{-1}{\Bbb E} [\langle f, \mu^{\xi}_{\la}
  \rangle]
  \to \tau E(\tau) \int_{Q_1} f(x) dx \ee
  where
  \be \label{thermolimit}  E(\tau):= E^{\xi}(\tau):=
  {\Bbb E} \left[\xi({\bf 0},{\cal P}^{\Psi}) \exp(-\Delta({\bf 0},{\cal P}^{\Psi}))\right].
  \ee
Moreover, if (\ref{mom}) is satisfied for some
$p>2$    
then
$\la^{-1} \langle f,\mu^{\xi}_{\la} \rangle$ converges to $\tau
E(\tau) \int_{Q_1} f(x) dx$ in $L^2.$
\end{theo}

 Note that $E(\tau)$ depends on the underlying intensity $\tau$ via
 ${\cal P}^{\Psi}$ even though this parameter does not explicitly show up
 in the defining formula. 
 Before stating variance asymptotics write
 $$\sigma^{\xi}[{\bf 0}] := {\Bbb E}\left[\xi^2({\bf 0},{\cal P}^{\Psi}) \exp(-\Delta({\bf 0},{\cal P}^{\Psi}))\right]$$
 and for all $x \in \R^d$ define the two
 point correlation functions for the functional $\xi$ over
 the Gibbsian input ${\cal P}^{\Psi}$ by
 \begin{equation}\label{SigmaDef}
  \sigma^{\xi}[{\bf 0},x] := {\Bbb E}\left[\xi({\bf 0},{\cal P}^{\Psi} \cup \{ x \})
   \xi(x,{\cal P}^{\Psi} \cup \{ {\bf 0} \}) \exp(-\Delta(\{ {\bf 0}, x \}, {\cal
   P}^{\Psi}))\right]
   - [E^{\xi}(\tau)]^2,
\end{equation}
where $\Delta(\{x,y\},\X) := \Psi(\X \cup \{ x, y \}) - \Psi(\X).$

\begin{theo}{(Variance asymptotics)}\label{VAR}
 Assume that $\xi$ is exponentially stabilizing and satisfies the $p$-moment
 condition (\ref{mom}) for some $p > 2$.  We have
 for each $f \in \B(Q_1)$
\be \label{varlimit}  \lim_{\la \to \infty} \la^{-1} \Var[\langle
f, \mu^{\xi}_{\la} \rangle] =  \tau V^{\xi}(\tau) \int_{Q_1}
f^2(x) dx, \ee
 where
 \be \label{varlimitdef} V^{\xi}(\tau) :=  \sigma^{\xi}[{\bf 0}]
 + \tau \int_{{\Bbb R}^d} \sigma^{\xi}[{\bf 0},x] dx < \infty. \ee
\end{theo}



Letting  $N(0, \sigma^2)$ denote a mean zero normal random
variable with variance $\sigma^2$, we have:

\begin{theo}{(CLT)}\label{CLT}
 Assume that $\xi$ is exponentially stabilizing and satisfies the $p$-moment
 condition (\ref{mom})
 for some $p > 2$. 
 We have for each $f \in \B(Q_1)$
 \be \label{1dimCLT}  \la^{-1/2}\langle f, \bar{\mu}^{\xi}_{\la} \rangle \tod N\left(0, \tau V^{\xi}(\tau) \int_{Q_1} f^2(x) dx\right),
 \ee
and the finite-dimensional distributions $( \la^{-1/2}\langle f_1,
\bar{\mu}^{\xi}_{\la} \rangle,...,\la^{-1/2}\langle f_m,
\bar{\mu}^{\xi}_{\la} \rangle),  \ f_1,...,f_m \in \B(Q_1),$
converge to those of a mean zero Gaussian field with covariance
kernel \be\label{covarr} (f_1, f_2) \mapsto \tau V^{\xi}(\tau)
\int_{Q_1} f_1(x)f_2(x)dx. \ee Moreover, if (\ref{mom}) is
satisfied for some $p > 3$ 
 then for all $\la \geq 2$ and all $f \in \B(Q_1)$
we have \be \label{rates} \sup_{t \in \R} \left| {\Bbb P} \left[
{\langle f, \bar{\mu}^{\xi}_{\la} \rangle \over \Var [\langle f,
\bar{\mu}^{\xi}_{\la} \rangle] } \leq t \right] - {\Bbb P} [N(0,1)
\leq t] \right| \leq C (\log \la)^{3d} \la^{-1/2}. \ee
\end{theo}

\ \

Assuming  that $\hat\xi(\cdot,\cdot;\la)$ is  a bounded
perturbation of a stabilizing functional $\xi$ we have:

 \begin{theo}\label{BDDLIMIT} The
analogs of Theorems \ref{WLLN}, \ref{VAR} and \ref{CLT} hold
  for $\xi$ replaced by $\hat\xi(\cdot,\cdot;\la)$ under their
  respective stabilization and moment assumptions.
 \end{theo}

\ \




{\bf Remarks.} {\em (i) Comparison with  \cite{FFG3}}.  The
results of \cite{FFG3} establish limit theory for functionals
$\xi$ of weakly dependent Gibbsian input, but essentially these
results require $\xi$ to
have finite range (finite range test functions). 
 Theorems \ref{WLLN}-\ref{BDDLIMIT}
extend \cite{FFG3} to cases when $\xi$ has infinite range and
stabilizes.

{\em (ii) Comparison with functionals on Poisson input.} Theorems
\ref{WLLN}-\ref{BDDLIMIT} show that the established limit theory
for stabilizing functionals on homogeneous Poisson input
\cite{BY2,Pe1, PY1,PY2,PY4} is insensitive to weak Gibbsian
modifications of the input.
 Thus the entirety of weak laws of large numbers and central limit theorems
 for functionals defined on homogeneous Poisson input
 given previously in the literature \cite{BY2}, \cite{Pe1}-\cite{PY4} extend
 to the corresponding analogous
results for functionals on point processes whose local
specification  (\ref{GIBBS1}) with respect to the Poisson process
is exponentially localized. 
If the input ${\cal P}^{\Psi}$ is Poisson, then the term
$\Delta(x,{\cal P}^{\Psi})$ vanishes, and hence  Theorem
\ref{WLLN} extends the Poisson weak law of large numbers in
Theorem 2.1 of \cite{PY4}. Likewise, Theorem \ref{VAR} extends the
variance asympotics of \cite{BY2} and \cite{Pe1}, whereas Theorem
\ref{CLT} extends the central limit theory of \cite{BY2},
\cite{Pe1} and  \cite{PY5}.

{\em (iii) Numerical evaluation of limits}. We emphasize that the
point process ${\cal P}^{\Psi}$ is intrinsically algorithmic; this
algorithmic scheme provides an exact (perfect) sampler
\cite{FFG3}.
It is computationally  efficient and yields a  numerical
evaluation of the limits (\ref{thermolimit}) and (\ref{varlimit}).

{\em (iv) Extensions and generalizations.}  The variance
convergence (\ref{varlimit}) and the asymptotic normality
(\ref{1dimCLT}) hold under weaker stabilization assumptions such
as power-law stabilization (see Penrose \cite{Pe1}), but the
resulting additional technical details obscure the the main ideas
of our approach, and thus we have not tried for the weakest
possible stabilization conditions.
 Similarly, counterparts to Theorems \ref{WLLN}-\ref{BDDLIMIT} should hold
for functionals
defined in terms of non-homogenous Gibbsian input, but we do not
provide the technical details here either.

\section{Proofs of main results}
This section is organized as follows.
 First, in Subsection \ref{PROP} we establish
 exponential clustering properties for stabilizing functionals of processes with
 exponentially localized potentials. 
Exponential clustering is central to our approach, as it shows
that the cumulants of $\langle f, \bar{\mu}_\la^{\xi} \rangle, f
\in \B(Q_1),$ converge to those of a normal random variable, that is to say they
vanish asymptotically upon suitable re-scaling for all orders above two.
 Then, in Subsections \ref{PROOFS1}
 and \ref{PROOFS2} we establish Theorems \ref{WLLN}, \ref{VAR} and Theorem \ref{CLT}
 respectively, using either the cumulant techniques developed in
 \cite{BY2} or the Stein techniques of \cite{PY5}. Subsection
 \ref{BDDLIMITsec} provides the proof of Theorem \ref{BDDLIMIT}.

\subsection{Exponential clustering lemma}\label{PROP}
 Let ${\cal P}^{\Psi}$ be a point process with exponentially localized potential and
 assume that $\xi$ is an exponentially stabilizing functional in the wide sense.
 For each Poisson-like configuration $\Xi$ we denote by $\xi[\Xi]$ this point
 configuration marked with the values of $\xi,$ that is to say each $x \in \Xi$
 carries the mark $\xi(x,\Xi).$ We have then
 \begin{lemm}\label{CLUSTERING}
  For each $k \geq 1$ there exist $M > 0$ and $c:= c(k) > 0$ such that for any
  deterministic points $x_1,\ldots,x_k \in \R^d$ the total variation distance between
  $\xi[{\cal P}^{\Psi}]$ restricted to the union $B_1(x_1) \cup \ldots \cup B_1(x_k)$
  and the product of respective restrictions of $\xi[{\cal P}^{\Psi}]$ to
  $B_1(x_1),\ldots,B_1(x_k)$ does not exceed $M k \exp(-c \min_{i,j} \dist(x_i,x_j)).$
 \end{lemm}

 {\em Proof.} The statement of the lemma is a consequence of the
  graphical construction of the process ${\cal P}^{\Psi}$ and of the exponential
  stabilization of $\xi.$ To see it, observe that the considered total variation
  distance does not exceed the probability of the event that the random sets
  $$ A_i := An_0[\bigcup_{x \in {\cal P}^{\Psi} \cap B_1(x_i)} B_R(x)]$$
  are not all disjoint, where $R:= R[x,{\cal P}^{\Psi}]$ and where
  the causal ancestor cone $An_0[A]$ is defined in section
  \ref{graphconstruct}. Indeed, if all $A_i$'s are
  disjoint then the values of $\xi$
  over all points in balls $B_1(x_i)$ depend on disjoint and hence independent portions
  of the free  birth-and-death process in the graphical construction. To
  complete the proof it suffices now to show
  that the probability ${\Bbb P}[A_i \cap A_j \neq \emptyset]$ decays exponentially
  with the distance between $x_i$ and $x_j$ for each $i$ and $j.$ Now, this follows because
  \begin{itemize}
   \item The number of points in $\P^{\Psi} \cap B_1(x_i)$ and $\P^{\Psi} \cap B_1(x_j)$
         admits super-exponentially decaying tails in view of the Poisson domination
         property of the \sr \ process $\P^{\Psi},$
   \item For each such point $x$ the stabilization radius $R[x,\P^{\Psi}]$ admits
         exponentially decaying tails by the wide sense exponential stabilization
         (\ref{ExpStab}),
   \item Consequently, the diameter of the union $\bigcup_{x \in {\cal P}^{\Psi} \cap B_1(x_i)} B_R(x)$
         of such balls has exponentially decaying tails too,
   \item Finally, using the exponential decay relation (\ref{ExpDec2}) for
         causal ancestor clan sizes in the graphical construction, we conclude
         that the diameter of $A_i$ also has exponentially decaying diameter.
  \end{itemize}
  The proof is hence complete. \qed

\subsection{Proof of Theorems \ref{WLLN} and \ref{VAR}}\label{PROOFS1}


There are several ways to prove
limit theorems for stabilizing translation invariant functionals.
To illustrate the new features arising in the setting of
functionals of Gibbsian input, we will first assume that $f$ is
a.e. continuous, that $\xi$ satisfies the moment condition
(\ref{mom}) for $p = 4$, and appeal to cumulant methods. In this
setting we may directly apply the cumulant methods developed in
Section 4
 of \cite{BY2}  (especially those methods used for proving statements (i) and (ii) of
 Theorem 2.1 there)  and hence
 we only provide crucial points, referring the reader to \cite{BY2}
 for further details. The arguments in Section 4 there show that
 \begin{equation}\label{LimRel1}
  \lim_{\la\to\infty} \la^{-1} {\Bbb E}\langle f, \mu^{\xi}_{\la} \rangle
    = c(\0) \int_{Q_1} f(u) du
 \end{equation}
 and
 \begin{equation}\label{LimRel2}
  \lim_{\la\to\infty} \la^{-1} \Var[\langle f, \mu^{\xi}_{\la} \rangle] =
    [q(\0) + \int_{\R^d} c(\0,x) dx] \int_{Q_1} f^2(u) du,
 \end{equation}
 where the correlation functions $c(x), q(x)$ and $c(x,y),\;x,y \in
 \R^d$,
 are the respective Radon-Nikodym derivatives given
 by
 $$ {\Bbb E}\mu^{\xi}(dx) = c(x)dx,$$ $$\E[(\mu^{\xi}(dx))^2] = q(x)dx \mbox{ and }$$
 $$ \Cov[\mu^{\xi}(dx),\mu^{\xi}(dy)] = c(x,y) dx dy,\; x \neq y, $$
 where $\mu^{\xi}$ is the infinite-volume empirical measure defined in (\ref{MUX}).
 Indeed,
 the main idea as briefly sketched below is to use stabilization,
 under the guise of
 the exponential clustering Lemma \ref{CLUSTERING} here, to show that
 when proving our results, in the $\la \to \infty$ limit we can
 safely
 replace (modulo a correction of order $o(\la^{1/2})$)
 the considered expression $\langle f, \mu^{\xi}_{\la} \rangle$
 by $\langle f, [\mu_{\infty}]^{\xi}_{\la} \rangle$ where
 $[\mu_{\infty}]^{\xi}_{\la} := \sum_{x \in {\cal P}^{\Psi} \cap Q_{\la}}
  \xi(x,{\cal P}^{\Psi}) \delta_{x/\la^{1/d}}.$ Now, the last expression
 coincides with $\int_{Q_{\la}} f_{\la} d \mu^{\xi},$ where $f_{\la}(x)
 := f(\la^{-1/d}x),\; x \in Q_{\la}.$ Consequently,
 ${\Bbb E}\langle f, \mu^{\xi}_{\la} \rangle =
  \langle f, {\Bbb E}\mu^{\xi}_{\la} \rangle$ for large $\la$ is well
 approximated by $\int_{Q_{\la}} f_{\la} d{\Bbb E}\mu^{\xi} =
 \int_{Q_{\la}} f_{\la}(x) c(x) dx = \la c(\0) \int_{Q_1} f(u) du$
 by translation invariance of $\mu^{\xi}$ and upon a variable change
 $u := \la^{-1/d} x.$ This yields (\ref{LimRel1}) for continuous
 $f.$
 To get (\ref{LimRel1}) in the general set-up, that is to say
 when $f \in \B(Q_1)$ and when $\xi$ satisfies the bounded moment condition (\ref{mom})
 for some $p > 1,$ one can follow the approach of \cite{Pe1}.
 Likewise, under exponential stabilization, $\Var[\langle f, \mu^{\xi}_{\la} \rangle]$
 is well approximated by $\Var[\int_{Q_{\la}} f_{\la} d \mu^{\xi}],$
  which by Campbell's theorem, equals
 $\int_{Q_{\la} \times Q_{\la}} f_{\la}(x) f_{\la}(y) \Var \mu^{\xi}(dx,dy),$
 where $\Var \mu^{\xi} := {\Bbb E} [\mu^{\xi} \otimes \mu^{\xi}] -
 [{\Bbb E}\mu^{\xi}]\otimes [{\Bbb E}\mu^{\xi}]$ is the variance
 measure of $\mu^{\xi};$ see Section 4 in \cite{BY2} and references
 therein for more details on moment measures. Using the usual decomposition
 of the variance measure into the diagonal and off-diagonal component
 \cite{BY2}, we see that the last expression equals
 \be \label{varex} \int_{Q_{\la}} f^2_{\la}(x) q(x) dx + \int_{Q_{\la} \times Q_{\la}}
  f_{\la}(x) f_{\la}(y) c(x,y) dx dy.\ee
  Using the continuity of $f_\la$, the translation invariance of $c(x,y)$ and
 $q(x)$, and  the exponential decay of $c(x,y)$ in the distance
 between $x$ and $y$ as guaranteed by the exponential clustering
 Lemma \ref{CLUSTERING} with $k = 2$, we come to (\ref{LimRel2}) as required.
 To show (\ref{LimRel2}) when $f \in \B(Q_1)$ and when $\xi$
 satisfies the moment condition (\ref{mom}) for some $p > 2$, we may modify the
 approach of \cite{Pe1}.

 Now, to calculate the correlation functions $c(\cdot), c(\cdot,\cdot)$ and
 $q(\cdot)$ in (\ref{LimRel1}) and
 (\ref{LimRel2}), note first that, given ${\cal P}^{\Psi}$
 in $\R^d \setminus dx,$ the probability of
 observing an extra point of ${\cal P}^{\Psi}$ at $x$ is
 $\tau \exp(-\Delta(x,{\cal P}^{\Psi})) dx$ as determined by the construction of
 the process in Subsection \ref{graphconstruct}, where $\tau dx$ corresponds to
 the birth attempt intensity at $x$ whereas $\exp(-\Delta(x,{\cal P}^{\Psi}))$
 comes from the acceptance probability. Consequently,
 ${\Bbb E}\mu^{\xi}(dx) = \tau {\Bbb E}\xi(x,{\cal P}^{\Psi})
 \exp(-\Delta(x,{\cal P}^{\Psi})) dx$ and hence
 \begin{equation}\label{CorrExp1}
  c(x) = \tau {\Bbb E}\xi(x,{\cal P}^{\Psi}) \exp(-\Delta(x,{\cal P}^{\Psi})).
 \end{equation}
 Likewise,
 \begin{equation}\label{CorrExp2}
  q(x) = \tau {\Bbb E}\xi^2(x,{\cal P}^{\Psi}) \exp(-\Delta(x,{\cal P}^{\Psi})).
 \end{equation}
 Further, for $x,y \in {\Bbb R}^d,$ given ${\cal P}^{\Psi}$ in $\R^d \setminus
 (dx \cup dy),$ the probability of observing extra points of ${\cal P}^{\Psi}$
 at $x$ and $y$ respectively is $\tau^2 \exp(-\Delta(\{x,y\},{\cal P}^{\Psi})) dx dy,$
 where again $\tau dx$ and $\tau dy$ are the probabilities that the birth attempts
 at $x$ and $y$ were made whereas $\exp(-\Delta(\{x,y\},{\cal P}^{\Psi}))$ is
 the probability that they were both accepted. Consequently,
 \begin{equation}\label{CorrExp3}
  c(x,y) = \tau^2 {\Bbb E}\xi(x,{\cal P}^{\Psi} \cup \{ y \})\xi(y,{\cal P}^{\Psi} \cup \{ x\})
           \exp(-\Delta(\{x,y\},{\cal P}^{\Psi})) - c(x) c(y).
 \end{equation}
 In other words, $c(x,y) = \tau^2 \sigma^{\xi}[{\bf 0},y-x]$ with
 $\sigma^{\xi}[\cdot,\cdot]$ as in (\ref{SigmaDef}).  
 The required relations (\ref{wllnlimit}) and (\ref{varlimit}) follow now by
 putting (\ref{LimRel1}) and (\ref{LimRel2}) together with (\ref{CorrExp1}),
 (\ref{CorrExp2}) and (\ref{CorrExp3}) and comparing with (\ref{thermolimit})
 and (\ref{varlimit}). The exponential clustering Lemma \ref{PROP} when
 combined with the moment conditions imposed on $\xi$ implies that the two-point correlation
 $c(x,y)$ exhibits exponential decay in the distance between $x$ and $y$ whence the
 integral $\int_{\R^d} c(\0,x) dx$ is finite. This observation
 allows us to conclude the finiteness of  $E(\tau)$ and $V(\tau),$ as  given by
  (\ref{thermolimit}) and (\ref{varlimitdef}), respectively. Consequently, the $L^2$-convergence
 stated in Theorem \ref{WLLN} follows now by the variance convergence in Theorem
 \ref{VAR} and, given (\ref{LimRel1}) and (\ref{LimRel2}),
 the proof of both of these theorems is complete. \qed

\subsection{Proof of Theorem \ref{CLT}}\label{PROOFS2}

When $f$ is continuous on $Q_1$ and when $\xi$ satisfies the
moment condition (\ref{mom}) for all $p$, the exponential
clustering Lemma \ref{CLUSTERING} allows us to use the techniques
 developed in Section 5 of \cite{BY2}, where it replaces the clustering Lemma 5.2,
 to show that all cumulants of $\langle f, \bar\mu^{\xi}_{\la} \rangle$ are all
 of the volume order $\la$ and hence, upon the $\la^{-k/2}$-re-scaling with $k$
 being the order of the cumulant, the cumulants of order higher than two vanish
 asymptotically yielding the required Gaussian limit; see \cite{BY2} for details.

 More generally,  for $f \in \B(Q_1)$ and when  $\xi$ satisfies the moment condition (\ref{mom})
for all $p > 2$,  the rate  (\ref{rates}) holds by following {\em
verbatim} the
 the Stein approach of \cite{PY5}, using wide sense stabilization
 and the exponential clustering Lemma \ref{CLUSTERING} instead of
 stabilization. Combining  (\ref{varlimit})  and (\ref{rates}) yields (\ref{1dimCLT}) for $f \in \B(Q_1)$.
This completes the proof of Theorem \ref{CLT}. \qed

\subsection{Proof of Theorem \ref{BDDLIMIT}} \label{BDDLIMITsec}
 The uniformly decaying bound on all moments of the perturbation
 term $\delta(\cdot,\cdot;\la)$ in (\ref{BDDIST1}), as stated
 in (\ref{BDDIST2}), combined 
 with the exponential stabilization
 of the perturbation, allows us to use
 the H\"older and Minkowski inequalities to conclude that the addition of $\delta(\cdot,\cdot;\la)$
 to $\xi$ does not affect the asymptotic behavior of the first and second order correlation
 functions. This means that the cumulant-based argument for Theorems \ref{WLLN} and \ref{VAR}
 carries over also for $\hat\xi(\cdot,\cdot;\la)$ with no further modifications. This yields
 the bounded perturbed versions of Theorems \ref{WLLN} and \ref{VAR} for continuous test
 functions and under the moment condition (\ref{mom}) with $p = 4$. To relax the moment
 conditions as in the respective statements and to get the results for general bounded
 test functions we resort again to the approach of \cite{Pe1}, which completes the proof
 of these two theorems for functionals with bounded perturbations. The remaining CLT Theorem
 \ref{CLT} for functionals with bounded perturbations follows now easily by the stabilization
 property imposed on the perturbation term $\delta(\cdot,\cdot;\la)$ in full analogy with
 the respective proof of Theorem \ref{CLT}.
 \qed

 \section{ Examples of Gibbs point processes with exponentially localized potentials}\label{S41}
\allco



 The notion of an exponentially localized potential $\Psi$ is general and
 includes the following non-exhaustive list of the corresponding point
 processes ${\cal P}^{\Psi}.$
 If an energy functional $\Psi$ has finite interaction range  so
 that its add-one potential satisfies $\Delta(x, \X) = \Delta(x,
 \X \cap B_r(x))$ for some $r$, as would be the case in many
 examples considered by \cite{FFG3}, then clearly (\ref{ExpDec})
 is satisfied and there usually are natural ways of ensuring that
 the constant $C_1$ is large enough so that exponential localization
 and (\ref{ExpDec2}) hold as well. These include decreasing the
 intensity $\tau$ of the underlying Poisson process ${\cal P}$
 which corresponds to increasing $\Delta(\cdot,\cdot)$ by a
 positive constant (low reference intensity/density regime) as well as multiplying
 $\Psi$ and hence also $\Delta(\cdot,\cdot)$ by some small enough
 $\beta > 0$ (high temperature regime). The following list is not
 limited to finite range energy functionals.

 (i) {\em Strauss processes.}  A Strauss process involves perturbing a Poisson process
 according to an exponential of the number of pairs of points  closer than a fixed
 cutoff. For such processes the add-one potential depends only on points within
 the cut-off range and so $\psi(r)$ vanishes when $r$ exceeds this cut-off.

 (ii) {\em Point processes with pair potential function.} A large
 class of Gibbs point processes \cite{SKM}, known as pairwise
 interaction point processes and including the Strauss process,
 has Hamiltonian $\Psi(\X):= \sum
 \sum_{i < j} \phi(||x_i - x_j ||), \ \X := \{x_i\},$ with $\phi$
 bounded below, usually assumed to be positive by absorbing the
 offending constant into the intensity of the underlying Poisson
 process.
If the pair potential function $\phi$ has finite range, as would
be the case with the Strauss process, then the potential $\Psi$ is
localized since $\psi$ vanishes beyond the interaction range. On
the other hand, suppose the pair potential function has infinite
range, but satisfies the following strengthened superstability
condition: $\phi$ decays exponentially fast and $\phi(s) = +
\infty$ for $s \leq r_0,$ that is there is a hard-core exclusion
condition forbidding the presence of two points within distance
less than $r_0$, 
\cite{Rue70}. In this context then the point process ${\cal
P}^{\Psi}$ is easily verified to be exponentially localized as
soon as the intensity $\tau$  is low enough (low density regime)
or $\phi$ admits a sufficiently small upper bound on its
oscillations (high temperature regime).

(iii) {\em Area interaction point processes}.
  This is a germ grain
 process, where the grain shape is a fixed compact convex set and
 where the potential 
 at each Poisson germ is determined by a function of the
 intersection of the grains at that germ. As a special and simple instance,
 suppose that the energy functional $\Psi_D(\X)$ is a scalar multiple
 $\gamma$ of the volume of the union of the radius $r$ balls
 centered at points $x \in \X \cap D$. Then, for $\gamma$ small enough,
 the resulting area
 interaction process (consisting of `ordered' points for negative
 interaction parameter
 $\gamma$ and `clustered' points for positive interaction parameter $\gamma$) is
 exponentially localized. More general energy functionals involve an additive term
 representing a scalar multiple of the total number of points \cite{BL}, which
 can be alternatively absorbed into the intensity.
 As noted in \cite{BL}, area interaction processes
 plausibly model certain biological processes, including those
 where the realization of the process represents spatial locations of plants (or
 animals) consuming food within distance $r$. The energy
 functional is then a scalar multiple of the area of the food supplying region.
 These are described more
 fully on p. 9 of \cite{FFG3} and in \cite{BL}.

(iv) {\em Point processes defined by the continuum Widom-Rowlinson
model.} Another example of the point process ${\cal P}^{\Psi}$ is
that defined in terms of the continuum Widom-Rowlinson model from
statistical physics, see \cite{WRO} as well as \cite{GHM}. Here we
have fixed radii (say radius equal to $a$) spheres of two types,
say $A$ and $B$, with interpenetrating spheres of similar types
but hard-core exclusion between the two types. This defines a
point process whose potential is exponentially localized as soon
as the reference intensity is low enough, 
since the function $\psi(r)$ vanishes when $r > 2a$. It is known,
see ibidem, that the continuum Widom-Rowlinson admits an
equivalent reformulation in terms of single-species gas of
interpenetrating spheres which is area-interacting in the sense of
point (iii) above \--- this is seen by integrating out the
positions of $B$ particles and keeping track of the locations of
$A$-particles only. Likewise, upon forgetting the marks carried by
the particles in the two-species representation one gets the
so-called random cluster representation for the Widom-Rowlinson
model, see \cite{CCK} and \cite{GHM}, from which the law of the
Widom-Rowlinson model can be recovered by assigning independent
and equiprobable $A$ and $B$-labels to maximal connected clusters
of particles, whence the
name {\it random cluster model}. 
Theorems \ref{WLLN}-\ref{BDDLIMIT} are valid for all of these
equivalent models, provided the intensity is low enough, as
discussed above.


(v) {\em Point processes given by hard-core model.}
 An important and natural model falling into the framework of our
 theory is the so-called hard-core model with low enough reference intensity.
 In its basic version the hard-core model, extensively studied in
 statistical mechanics, arises by conditioning a Poisson point
 process on containing no two points within distance less than
 $2r$ for some $r>0$ standing for a parameter of the model.
 Clearly, this process admits a Gibbsian description with
 $\Psi$ set to $+\infty$ if there are two points closer than
 $2r$ from each other and $0$ otherwise. Consequently, the
 potential is exponentially localized if the reference intensity of
 the underlying Poisson point process is low enough or
 if $r$ is small enough, which also reduces to decreasing
 the reference intensity  upon appropriate re-scaling (in fact rather
 than imposing separate conditions on $r$ and the reference intensity
 $\tau$ it is enough to require that $\tau r^d$ be small
 enough, as easily checked by re-scaling).

 (vi) {\em Truncated Poisson process.} The hard-core gas is a
 particular example of a truncated Poisson process.
 In general, a truncated Poisson process arises by conditioning
 a Poisson point process on the event that a certain family of
 constraints is fulfilled. In this paper the constraints imposed
 are of the following form: we fix a certain family of bounded sets
 and require that none of these sets contain more than a certain
 given number of points. Such processes are used in modelling of
 communication networks \cite{BBK}.
 In particular, if we require that  no ball of radius $r$
 contain more than some constant number $k$ of Poisson points, then
 $\psi$ vanishes beyond $r$ and the associated point process has
 an exponentially localized potential, possibly upon decreasing
 the intensity.


\section{Applications} \label{applic}


Below we indicate some  applications of our main results. This
list is not exhaustive and does not include applications to e.g.
germ-grain models  where the germs arise as the realization of the
Gibbsian point process $\P^{\Psi}$ with an exponentially localized
potential.

\subsection{RSA packing with Gibbsian input}

Let $\X \subset \R^d$ be locally finite.  Consider a sequence of
unit volume $d$-dimensional Euclidean balls $B_{1}, B_{2},....$
with centers arriving sequentially at points in $\X$. The first
ball $B_{1}$ to arrive is {\em packed} and recursively, for $i
=2,3,...$ let the $ith$ ball be packed if it does not overlap any
ball in $B_{1}, B_{2},....B_{i-1}$ which has already been packed.
Let $\xi(x,\X)$ be either $0$ or $1$, depending on whether the
ball arriving at $x$ is either packed or discarded.

When $\X$ is the realization of a Poisson point process  on
$Q_\la$, this packing process is known as random sequential
adsorption (RSA) with Poisson input on $Q_\la$. When $\X$ is the
realization of an infinite sequence of independent random
$d$-vectors uniformly distributed on the cube $Q_\la$, then this
is called the RSA process with infinite binomial input; in such
cases,  RSA packing terminates when it is no longer possible to
pack additional balls. In dimension $d = 1$, this process is known
as the R\'enyi car parking problem \cite{Re}.
 In the infinite input setting and when $d = 1$  R\'enyi \cite{Re} (respectively
Dvoretzky and Robbins \cite{DR}) proved that the total number of
parked cars satisfies a weak law of large numbers (respectively
central limit theorem) as $\la \to \infty$; recently these results
were shown to hold for all dimensions in \cite{Pe1} and
\cite{SPY}. 

Virtually all limit results for RSA packing assume that the input
is either Poisson or a fixed number of independent identically
distributed random variables. To the best of our knowledge, RSA
packing problems with Gibbsian input have not been considered
before in the literature.  The following theorem widens
 the scope of the existing limit results for RSA
packing.  Put $$\mu_\la^{\xi}:= \sum_{x \in {\cal P}^{\Psi} \cap
Q_{\la}} \xi(x,{\cal P}^{\Psi} \cap
Q_{\la})\delta_{x/\la^{1/d}},$$ so that $N(\P^{\Psi} \cap Q_\la)
:= \sum_{x \in \P^{\Psi} \cap Q_\la} \xi(x, \P^{\Psi} \cap
 Q_\la)$  denotes the total number of balls packed on $Q_\la$ from
the collection of balls with centers in $\P^{\Psi} \cap Q_\la$.


\begin{theo} \label{packing} Let $\P^{\Psi}$ be Gibbsian input with an exponentially
localized potential. Then
 \be \label{pack1} \la^{-1} N(\P^{\Psi} \cap
Q_\la) \to \tau {\Bbb E} \left[\xi({\bf 0},{\cal P}^{\Psi})
\exp(-\Delta({\bf 0},{\cal P}^{\Psi}))\right] \ \ \ \text{in} \
L^2 \ee and
 \be \label{pack2} \la^{-1} \Var[N(\P^{\Psi} \cap Q_\la)] \to \tau
V^\xi(\tau) \ee where $V^\xi(\tau)$ is given by
(\ref{varlimitdef}). The finite-dimensional distributions $(
\la^{-1/2}\langle f_1, \bar{\mu}^{\xi}_{\la}
\rangle,...,\la^{-1/2}\langle f_m, \bar{\mu}^{\xi}_{\la} \rangle),
$ $f_1,...,f_m \in \B(Q_1),$ converge to those of a mean zero
Gaussian field with covariance kernel \be\label{pack3} (f_1, f_2)
\mapsto \tau V^{\xi}(\tau) \int_{Q_1} f_1(x)f_2(x)dx. \ee
\end{theo}

{\bf Remark.} As spelled out in \cite{PY2}, Theorem \ref{packing}
also applies to  related packing models, including spatial birth
growth models with Gibbsian input as well as RSA models with balls
replaced by particles of random size/shape/charge, and ballistic
deposition models.

 {\em Proof.}   The approach used in \cite{PY2} shows that  the
packing functional $\xi(x, \cdot)$ is exponentially stabilizing on
$\sr$ sets. Indeed, any Poisson-like set $\Xi$ can be coupled on
the common underlying probability space with a dominating Poisson
point process $\P_{\tau}$ of finite intensity $\tau$, $\tau$
large, and containing $\Xi$ a.s. Now, the idea underlying the
argument in \cite{PY2}  shows that the packing status of a point
$x$ in a configuration $\X$ depends on $\X$ only through its
algorithmically determined sub-configuration $Cl[x,\X]$ referred
to as the {\it causal cone} or {\it causal cluster} of $x$ in the
presence of $\X,$ see \cite{PY2} for details. The causal cluster
$Cl[x,\X]$ is easily seen to be non-decreasing in $\X.$ In
particular, using that $\Xi \subseteq \P_{\tau}$ yields $Cl[x,\Xi]
\subseteq Cl[x,\P_{\tau}]$ a.s. for $x \in \Xi.$
However, by the arguments in section 4 of \cite{PY2}, the causal
clusters generated by points of $\P_{\tau}$ exhibit exponential
decay, and hence so do causal clusters of points in $\Xi$ showing
that the packing functional $\xi(x, \cdot)$ is exponentially
stabilizing on Poisson-like sets, in particular on $\P^{\Psi}$. In
other words, $\xi(x,\cdot)$ is exponentially stabilizing in the
wide sense. Clearly $\xi$ satisfies the bounded moment condition
(\ref{mom}) and therefore Theorems \ref{WLLN}- \ref{CLT} show that
the $\langle f, \bar{\mu}^{\xi}_{\la} \rangle, \ f \in \B(Q_1)$,
satisfy the weak law of large numbers and central limit theorem
given by (\ref{pack1}- \ref{pack3}), respectively.  \qed

\subsection{Functionals of Euclidean graphs on Gibbsian input}\label{Euclid}

In many cases, showing exponential stabilization of functionals of
geometric graphs over Poisson point sets \cite{BY2,PY1},  can be
reduced to upper bounding the probability that 
 regions in $\R^d$ are devoid of points  by a term which
decays exponentially with the volume of the region. When the
underlying point set is Poisson, as in \cite{BY2,PY1}, then we
obtain the desired exponential decay.  When the underlying point
set is Poisson-like, the desired exponential decay is an immediate
consequence of condition (\ref{LBd}).   In this way  the existing
stabilization proofs for functionals over Poisson point sets carry
over to stabilizing functionals on  Poisson-like point sets. This
extends central limit theorems for functionals of Euclidean graphs
on Poisson input to the corresponding central limit theorems for
functionals defined over Gibbsian input.   The following
applications illustrate this.


\vskip.3cm

{\bf (i) $k$-nearest neighbors graph.} The $k$-nearest neighbors
(undirected) graph on the vertex set $\X$, denoted $NG(\X)$, is
defined to be the graph obtained by including $\{x,y\}$ as an edge
whenever $y$ is one of the $k$ nearest neighbors of $x$ and/or $x$
is one of the $k$ nearest neighbors of $y$.  The $k$-nearest
neighbors (directed) graph on $\X$, denoted $NG'(\X)$, is obtained
by placing a directed edge between each point and its $k$-nearest
neighbors.

{\em Total edge length of $k$-nearest neighbors graph.} Let
$L(\X)$ denote the total edge length of $NG(\X)$ and let
$\xi(x,\X)$ denote one half the sum of the edge lengths of edges
in $NG(\X)$ which are incident to $x$. Put
$$\mu^{\xi}_{\la} := \sum_{x \in {\cal P}^{\Psi} \cap Q_{\la}}
\xi(x,{\cal P}^{\Psi} \cap Q_{\la})\delta_{x/\la^{1/d}}.$$



\begin{theo} \label{knnedgelenththm} Let $\P^{\Psi}$ be Gibbsian input with exponentially
localized potential. Then
 \be \label{nned1} \la^{-1} L(\P^{\Psi} \cap
Q_\la) \to \tau {\Bbb E} \left[\xi({\bf 0},{\cal P}^{\Psi})
\exp(-\Delta({\bf 0},{\cal P}^{\Psi}))\right] \ \ \ \text{in} \
L^2 \ee and
 \be \label{nned2} \la^{-1} \Var[L(\P^{\Psi} \cap Q_\la)] \to \tau
V^\xi(\tau) \ee where $V^\xi(\tau)$ is given by
(\ref{varlimitdef}). The finite-dimensional distributions $(
\la^{-1/2}\langle f_1, \bar{\mu}^{\xi}_{\la}
\rangle,...,\la^{-1/2}\langle f_m, \bar{\mu}^{\xi}_{\la}
\rangle),$ $ f_1,...,f_m \in \B(Q_1),$ converge to those of a mean
zero Gaussian field with covariance kernel \be\label{nned3} (f_1,
f_2) \mapsto \tau V^{\xi}(\tau) \int_{Q_1} f_1(x)f_2(x)dx. \ee
\end{theo}

{\bf Remark.} 
Theorem \ref{knnedgelenththm} generalizes Theorem 6.1 of
\cite{PY1}, which is restricted to nearest neighbor graphs defined
on Poisson input.

\vskip.5cm

{\em Proof.}  Let $\Xi$ be a Poisson-like point set. Considering
the arguments in the proofs of Theorem 6.1 of \cite{PY1} and
Theorem 3.1 of \cite{BY2},
it is easily seen that the set of edges incident to any point $x$
in $NG(\Xi)$ is unaffected by the addition or removal of points
outside a ball of random radius $R$. Moreover, the radius $R$ has
exponentially decaying tails, which may be seen as follows. For
simplicity we prove exponential stabilization in dimension two,
but the argument is easily extended to higher dimensions by using
cones instead of triangles (for $d = 1$ we use intervals instead
of triangles).
 For each $t > 0$ construct $6(k + 1)$
triangles  $T_j(t), \ 1 \leq j \leq 6(k + 1)$, such that $x$ is a
vertex of each triangle and such that each triangle with edge
containing $x$ has length $t$. Let $R_x$ be the minimum $t$ such
that each triangle contains at least one point from $\Xi$. In such
a situation, the union of the $6(k + 1)$ triangles $T_j(t), \ 1
\leq j \leq 6(k + 1),$ may be partitioned into $6$ equilateral
triangles with common edge length $t$, each triangle containing at
least $k + 1$ points. Then, because $\Xi$ is $\sr$, it follows
that ${\Bbb P} [R_x \geq t] \leq 6(k + 1) \exp(-C_3t^d)$.
Moreover, as explained in the proof of Theorem 6.1 of \cite{PY1},
simple geometry shows that $4R_x$ is a radius of stabilization for
the functional $\xi$ at $x$. Thus $\xi$ is exponentially
stabilizing.

 An easy
modification of the proof of Lemma 6.2 of \cite{PY1} shows that
$L$ moreover satisfies the $p$-moments condition (\ref{mom}) for
all $p$. Therefore Theorems \ref{WLLN}- \ref{CLT} show that the
$\langle f, \bar{\mu}^{\xi}_{\la} \rangle, \ f \in \B(Q_1)$,
satisfy the weak law of large numbers and central limit theorem
given by (\ref{nned1}- \ref{nned3}), respectively.  \qed
\vskip.3cm



{\em Number of components in nearest neighbors graph. }
Let $k = 1$.  Given a locally finite point set $\X$, let
$\xi^{[c]}(x, \X)$ denote the reciprocal of the cardinality
of the component in $NG(\X)$ which contains $x$.  Thus $H(\X):= \sum_{x
\in \X} \xi^{[c]}(x, \X)$ denotes the total number of components
of $NG(\X)$. Put $$\mu_\la^{\xi}:= \sum_{x \in {\cal P}^{\Psi}
\cap Q_{\la}}  \xi^{[c]}(x,{\cal P}^{\Psi} \cap
Q_{\la})\delta_{x/\la^{1/d}}.$$

\begin{theo} \label{knnthm} Let $\P^{\Psi}$ be Gibbsian input with exponentially
localized potential. Then
 \be \label{nnc1} \la^{-1} H(\P^{\Psi} \cap
Q_\la) \to \tau {\Bbb E} \left[\xi({\bf 0},{\cal P}^{\Psi})
\exp(-\Delta({\bf 0},{\cal P}^{\Psi}))\right] \ \ \ \text{in} \
L^2 \ee and
 \be \label{nnc2} \la^{-1} \Var[H(\P^{\Psi} \cap Q_\la)] \to \tau
V^\xi(\tau) \ee where $V^\xi(\tau)$ is given by
(\ref{varlimitdef}). The finite-dimensional distributions $(
\la^{-1/2}\langle f_1, \bar{\mu}^{\xi}_{\la}
\rangle,...,\la^{-1/2}\langle f_m, \bar{\mu}^{\xi}_{\la} \rangle),$
$ f_1,...,f_m \in \B(Q_1),$ converge to those of a mean zero
Gaussian field with covariance kernel \be\label{nnc3} (f_1, f_2)
\mapsto \tau V^{\xi}(\tau) \int_{Q_1} f_1(x)f_2(x)dx. \ee
\end{theo}

{\em Proof.} We  establish that $\xi^{[c]}$ is exponentially
stabilizing on Poisson-like sets $\Xi$ and appeal to Theorems
\ref{WLLN}-\ref{CLT}. When $k = 1$, the \sr \ properties of the
input process and the methods of H$\ddot{a}$ggstr$\ddot{o}$m and
Meester \cite{HaM} and Kozakova, Meester, and Nanda \cite{KoMeNa}
show there are no infinite clusters in $NG(\Xi)$.
Moreover, the proof of Theorem 1.1, 1.2 and Proposition 2.2 of
\cite{KoMeNa} and property (\ref{LBd}) of \sr  \ processes
show that the finite clusters in
$NG(\Xi)$ have (super)exponentially decaying cardinalities and
diameters. Now we show exponential stabilization of $\xi$ as
follows. Let $R_x$ be the radius of the cluster in $NG(\Xi)$
containing $x$. Write $B_r$ for $B_r(\0)$. Put
$$
R:= \sup_{x \in B_{R_\0} \cap \Xi} R_x.
$$
It is not hard to see  that $R$ has
exponentially decaying tail.  Indeed, writing
$$
P[ R > t] = \sum_{j =1}^{\infty} P [\sup_{x \in B_{R_\0} \cap \Xi}
R_x > t, \ 2^{-j-1} \leq R_\0 < 2^j]
$$
$$
 \leq \sum_{j =1}^{\infty} P
[\sup_{x \in B_{2^j} \cap \Xi} R_x > t, \ 2^{-j-1} \leq R_\0 ]
\leq \sum_{j =1}^{\infty} \exp(-C 2^j) P [\sup_{x \in B_{2^j} \cap
\Xi} R_x > t]
$$
and then noting that the cardinality of ${x \in B_{2^j} \cap \Xi}$
decays polynomially fast in $2^j$ with overwhelming probability,
we obtain the desired exponential decay of $R$. We can also show
that $4R$ is a radius of stabilization for $\xi^{[c]}$ at the
origin (see proof of Lemma 6.1 of \cite{PY1}) . Since $\xi^{[c]}$
trivially satisfies the bounded moments condition (\ref{mom}) for
all $p$, the weak law and central limit theorem  for $H(\P^{\Psi}
\cap Q_\la)$  and $\la^{-1/2}\langle f_1, \bar{\mu}^{\xi}_{\la}
\rangle$ follows by Theorems \ref{WLLN}- \ref{CLT}.  \qed



\vskip.3cm


{\bf (ii) Voronoi tessellations}. Given $\X \subset \R^d$ and $x
\in \X$, the set of points in $\R^d$ closer to $x$ than to any
other point of $\X$ is a convex polyhedral cell $C(x,\X)$. The
collection of  cells $C(x,\X), x \in \X$, form a partition of
$\R^d$ which is termed the Voronoi tessellation induced by $\X$.

{\em Total edge length.} Given $\X \subset \R^2$, let $L(x, \X)$
denote one half the total edge length of the finite edges in the
cell $C(x,\X)$.  It is easy to see that $L$ is exponentially
stabilizing on Poisson-like sets $\Xi$. Indeed, when $d = 2$,  it
suffices to follow the arguments in the proof of Theorem 8.1 of
\cite{PY1} and to note that stabilization radius depends on
finding a minimum edge length $t$ such that 12 isosceles triangles
with this edge length have at least one point from $\Xi$ in them.
Because $\Xi$ is Poisson-like we may follow the arguments in
\cite{PY1} verbatim to see that $L$ stabilizes. See section 6.3 of
\cite{Pe1} for the case $d > 2$.

As in section 6.3 of \cite{Pe1}, we may show that $L$ also
satisfies the moment condition (\ref{mom}) for $p = 3$. It follows
that the total edge length
 $\sum_{x \in \P^{\Psi} \cap Q_\la} L(x, \P^{\Psi} \cap Q_\la)$
of the Voronoi tessellation on $\P^{\Psi} \cap Q_\la$ satisfies
the weak law and central limit theorem as $\la \to \infty$. In
other words, the Voronoi analog of Theorem \ref{knnedgelenththm}
holds.

\vskip.3cm

{\bf (iii) Other proximity graphs}. There are further examples
where showing exponential stabilization of functionals of
geometric graphs (in the wide sense)  involves  upper bounding the
probability that  regions in $\R^d$ are devoid of Poisson-like
points. Such estimates are available in the Poisson setting and it
is not difficult to extend them to Poisson-like point sets.

In this way, by modifying the methods of \cite{PY1} (sections 7
and 9) and \cite{BY2} (section 3.1), we obtain weak laws of large
numbers and central limit theorems for the total edge length of
the sphere of influence graph, the Delaunay graph, the Gabriel
graph, and the relative neighborhood graph over Gibbsian input
$\P^{\Psi}$.

\subsection{ Gibbsian continuum percolation}

Let $\X$ be a locally finite point set and connect all pairs of
points which are at most a unit distance apart. The resulting
graph is equivalent to the basic model of continuum percolation,
in which one considers the union of the radius $1$ balls centered
at points of $\X$, see Section 12.10 in \cite{Grimm}.
Let $\xi^{[c]}(x,\X)$ be the reciprocal of the
size of the component containing $x$, so that
$$
N(\X) = \sum_{x \in \X} \xi^{[c]}(x, \X)
$$
counts the number of finite components in  $G$.

Section 9 of \cite{PY1} discusses central limit theorems  for
$N(\P \cap Q_\la).$ 
Using Theorem  \ref{CLT} we can generalize these results to obtain
a central limit theorem for the number of
components $N(\P^{\Psi} \cap Q_\la)$ 
in the continuum percolation model on Gibbsian input in the
subcritical regime, possibly of interest in the context of sensor
networks on Gibbsian point sets. To formulate this central limit
theorem assume that our Gibbs point process $\P^{\Psi}$ with
exponentially localized potential is such that it admits a
stochastic upper bound by a homogeneous Poisson point process of
some intensity $\tau$ falling into the subcritical regime of the
considered continuum percolation (see Section 12.10 in
\cite{Grimm}). Note that due to the Poisson-like nature of the
input process $\P^{\Psi}$ this is always possible upon  spatial
re-scaling. We argue that $\xi^{[c]}$ is exponentially stabilizing
on Poisson-like sets $\Xi$ as follows. If $\Xi$ is Poisson-like
and if $\tau$ is subcritical, then $\Xi$ is also subcritical by
stochastic domination. Consequently, the diameter of the connected
cluster emanating from a given point has exponentially decaying
tails, see ibidem. This yields the required exponential
stabilization upon noting that $\xi^{[c]}(x,\cdot)$ does not
depend on point configurations outside the connected cluster at
$x.$
Moreover, $\xi^{[c]}$ is bounded above by one and thus satisfies
the moments condition (\ref{mom}). Hence by Theorems
\ref{WLLN}-\ref{CLT}, $N(\P^{\Psi} \cap Q_\la)$ satisfies the weak
law of large numbers and central limit theorem, exactly as in the
statement of Theorem \ref{knnthm}.

\subsection{Functionals on Gibbsian loss networks}

Given the Poisson point process $\P:= \P_\tau$, consider the
following Gibbs point process. Fix an integer $m \in \N$. Attach
to each point of $\P$ a bounded convex grain $K$ and put the
potential $\Psi$ to be infinite whenever the grain $K$ at one
point has non-empty intersection with more than $m$ other grains.
This condition prohibits overcrowding, and, for more general
repulsive models, one can put $\Psi$ large and finite whenever the
grain $K$ at one point has non-empty intersection with a large
number (some number less than $m$) of other grains. The resulting
point process, which we call $\P^{\Psi}$, represents a version of
spatial loss networks appearing in mobile and wireless
communications. As discussed in point (vi) in Subsection \ref{S41}
the so defined $\Psi$ is exponentially localized as soon as the
underlying intensity $\tau$ is small enough.


Let $\K$ be an open convex cone in $\R^d$ (a cone is a set that is
invariant under dilations) with apex at the origin.  Given $x,y
\in \P^{\Psi}$, we say that $y$ is {\em connected to $x$}, written
 $x \rightarrow y$, if there is a sequence of points
$\{x_i\}_{i=1}^n \in (\K + x) \cap \P^{\Psi}$, $|x_i - x_{i + 1}|
\leq 1$, $|x_1 - x| \leq 1$ and $ |y - x_{n + 1}| \leq 1$.  If the
length of this sequence does not exceed a given $m,$ we write $x
\rightarrow_m y.$ For all $r > 0$ let $B_r^{\K}(x):= x + (\K \cap
B_r(\0))$.

{\em Coverage functionals.} The functional  $\xi(x,\P^{\Psi}):=
\sup \{r \in \R: \  x \ \rightarrow y \ \ \text{for all} \ y \in
B_r^{\K}(x) \cap \P^{\Psi} \}$ determines the maximal coverage
range of the network at $x$ in the direction of the cone $\K$.

The coverage measure is $\mu^{\xi}_{\la}:= \sum_{x \in \P^{\Psi}
\cap Q_\la }  \xi (x, \P^{\Psi} \cap Q_\la ) \de_{x/\la^{1/d}} $.
Confining attention to $\P,$ where $\tau$ belongs to the the
subcritical regime for  continuum percolation,  $\P^{\Psi}$ is in
turn  subcritical  because of Poisson domination. Since the
continuum percolation clusters generated by any Poisson-like set
$\Xi$ have exponentially decaying diameter, it follows that $\xi$
stabilizes in the wide sense (recall the proof for the number of
components in the continuum percolation model) and that $\xi$
admits an
exponential moment.  
  By appealing to  Theorems
\ref{WLLN} and \ref{CLT}, we obtain a weak law of large numbers
and central limit theorem for both the coverage measure
$\mu^{\xi}_{\la}$ and the total coverage $\sum_{x \in \P^{\Psi}
\cap Q_\la } \xi(x, \P^{\Psi} \cap Q_\la )$.


{\em Network reach functional.} Say that the network has {\em
reach at least } $r$ at $x$ if $x \rightarrow y$ for all $y \in
B_r^{\K}(x) \cap \P^{\Psi} $. 
 Put $\xi_r(x,\P^{\Psi}):= 1$ if the network has
reach at least $r$ at $x$ and otherwise put $\xi_r(x,\P^{\Psi}):=
0$. Theorems \ref{WLLN} and \ref{CLT} yield a weak law of large
numbers and central limit theorem for the total network reach
$\sum_{x \in \P^{\Psi} \cap Q_\la } \xi_r(x, \P^{\Psi} \cap Q_\la
)$.

{\em Number of customers obtaining coverage.} Independently mark
each point $x$ of $\P^{\Psi}$ with mark  $T$ (transmitter) with
probability $p
> 0$ and with mark $R$ (receiver) with the complement probability.
Then define the reception functional $\xi(x, \P^{\Psi})$ to be $1$
if $x$ is marked with $T$ or (when $x$ is marked with $R$) if  $z
\rightarrow x$ for some $z$ in the transmitter set $\{ z \in
\P^{\Psi}:  \ z \ \text{marked with} \ T \}$. Put $\xi(x,
\P^{\Psi})$ to be zero otherwise.  Thus $\xi(x, \cdot ) $ counts
when a customer at $x$ gets coverage and the limit theory for the
sum $ \sum_{x \in \P^{\Psi} \cap Q_\la } \xi (x, \P^{\Psi} \cap
Q_\la )$, which counts the total number of receivers (customers)
obtaining network coverage, is given by Theorems \ref{WLLN} and
\ref{CLT}. 


{\em Connectivity functional.} Given a broadcast range $r
> 0$  and the transmitter
set $\{ z \in \P^{\Psi}:  \ z \ \text{marked with} \ T
\}$, let $\xi_r(x, \P^{\Psi})$ be the minimum number, say $m$,
such that every point in $y \in B_r^{\K}(x) \cap \P^{\Psi}$ can be
reached from some transmitter $z \in \P^{\Psi}$ with $m$ or fewer
edges or hops, that is to say there exists a transmitter $z$ such
that $z \rightarrow_m y$ for all $y \in B_r^{\K}(x) \cap \P^{\Psi}.$
 Thus all receivers in the broadcast range $r > 0$ can be linked to
 a transmitter in $m$ or fewer hops.  Small values of $\xi_r(x, \P^{\Psi})$
represent high network connectivity;  for each $r > 0$, Theorems
\ref{WLLN} and \ref{CLT} provide a weak law of large numbers and
central limit theorem for the connectivity functional $\sum
\xi_r(z,\P^{\Psi} \cap Q_\la )$ as $\la \to \infty$.

\section{Gibbsian quantization for non-singular probability measures}\label{QUANT}
\allco

 Quantization for probability measures concerns the best
approximation of a $d$-dimensional probability measure $P$ by a
discrete measure supported by a set $\X_n$ having $n$ atoms. It
involves a partitioning problem of the underlying space and it
arises in a variety of scientific fields, including information
theory, cluster analysis, stochastic processes, and mathematical
models in economics \cite{GL}. The goal is to optimally represent
$P$, here assumed non-singular with density $h,$ 
 with a point set $\X_n$, where
optimality involves minimizing the {\em $L^r$ stochastic
quantization error} 
 (or `random distortion  error') given by
$$
I(\X_n):= \int_{\R^d} (\min_{x \in \X_n} |y-x|)^r P(dy)=  \sum_{x
\in \X} \int_{C(x,\X_n)} |y - x|^r P(dy).
$$
Recall that for all $x$ and locally finite point sets $\X$,
$C(x,\X)$ denotes the Voronoi cell (`Voronoi quantizer') induced
by the Euclidean norm around $x$ with respect to $\X$.

 The optimal (non-random)
quantization error is given by $\min_{\X_n} I(\X_n)$ and the
seminal  work of Bucklew and Wise \cite{BW} shows that this
minimal error satisfies \be \label{BW} \lim_{n \to \infty} n^{r/d}
\min_{\X_n} I(\X_n) = Q_{r,d} ||h||_{d/(d+ r)}\ee  where
$||h||_{d/(d+ r)}$ denotes the $d/(d + r)$ norm of the density $h$
and where the so-called {\em $r$th
quantization coefficient} $Q_{r,d}$ is some positive constant not
known to have a closed form expression.

The first order asymptotics for the distortion error on i.i.d.
points sets (that is to say letting $\X_n$ consist of i.i.d.
random variables) was first investigated by Zador \cite{Za} and
later by Graf and Luschgy \cite{GL} and Cohort \cite{Co}. Letting
$X_n$ be i.i.d. random variables with common density
$h^{d/(d+r)}/\int h^{d/{d + r}}$ and $\omega_d$ the volume of the
unit radius $d$-dimensional ball in $\R^d$, Zador's theorem shows
\be \label{Za} \lim_{n \to \infty} n^{r/d} I(\X_n) =
\omega_d^{-r/d} \Gamma(1 + r/d) ||h||_{d/(d+ r)},\ee whence (see
Prop. 9.3 in \cite{GL}) the  upper bound \be \label{upbd} Q_{r,d}
\leq \omega_d^{-r/d} \Gamma(1 + r/d). \ee

 Molchanov
and Tontchev \cite{MT} have pointed out the desirability  for
quantization via Poisson point sets and our purpose here is to
establish asymptotics of the quantization error
 on Gibbsian input. This is done as follows.
For $\la > 0$ and a finite point configuration $\X$ we abbreviate
$\X_{(\la)} := \la^{-1/d} \X.$ Moreover, we write $\tilde\X := \X
\cap Q_1$ so that in particular $\tilde\X_{(\la)}:= \la^{-1/d} \X
\cap Q_1.$
Recall also that we write $\P^{\Psi}$ for
$\P^{\Psi}_{\tau}$ as in the previous sections. Consider the
random point measures
 induced by the distortion arising from
$\tilde\P^{\Psi}_{(\la)}$, namely
$$
\mu^{\Psi}_\la := \sum_{x \in \tilde\P^{\Psi}_{(\la)}}
 \int_{C(x,\tilde\P^{\Psi}_{(\la)}) } |y - x|^r P(dy) \delta_x. \ \
$$

%

We will be interested in the asymptotic behavior of the random
integrals $\langle f, \mu^{\Psi}_\la \rangle$.  Clearly, when $f
\equiv 1$ then $\langle f, \mu^{\Psi}_\la \rangle$ gives another
expression for the distortion $I(\tilde\P^{\Psi}_{(\la)})$.  On the other
hand,  if $f = 1_B$, then $\langle f, \cdot \rangle$ measures the
local
distortion. 
 This section establishes  mean
and variance asymptotics for $\langle f, \mu^{\Psi}_\la \rangle$
as well as convergence of the finite-dimensional distributions of
$\langle f, \mu^{\Psi}_\la \rangle$. Since we will no longer be
working with translation invariant $\xi$, we will need to appeal
to Theorem \ref{BDDLIMIT}.

Put
 $$ M^{\Psi}(\tau) := \int_{C(\0,\P^{\Psi})}
    |w|^r dw \exp(-\Delta({\bf 0},{\cal P}^{\Psi}))$$
 where, recall, $\tau$ is the intensity of the reference process
 $\P.$ Note that $M^{\Psi}(\tau)$ does depend on $\tau$ through
 $\P^{\Psi}.$
 Changing the order of integration we have
 \begin{equation}\label{EMcalc}
    \E M^{\Psi}(\tau) = \E \int_{\R^d} |y|^r {\bf 1}_{\P^{\Psi} \cap B_{|y|}(y) = \emptyset}
    \exp(-\Delta({\bf 0},\P^{\Psi})) dy =
 \end{equation}
 $$
    \int_{\R^d} |y|^r \E [\exp(-\Delta({\bf 0},\P^{\Psi}))
    {\bf 1}_{\P^{\Psi} \cap B_{|y|}(y) ) = \emptyset}] dy.
 $$
 In the special case where $\Psi \equiv 0$ (i.e. $\P^{\Psi}$ coincides with
 the reference process $\P$) and where the intensity $\tau$ of $\P$ is $1$
 we readily get $\E M^{0}(1) = \Gamma(1 + {r\over d}) \omega_d^{-r/d}.$
 More generally $\E M^0(\tau) = \tau^{-(1+r/d)}\Gamma(1 + {r\over d}) \omega_d^{-r/d}.$
Put
$$
V^{\Psi}(\tau) :=  \E [M^{\Psi}(\tau))^2] \ + $$
$$ \int_{\R^d } \left( \E \left[ \int_{C(\0, \P^{\Psi}
\cup \{ y \} )} |w|^r dw \int_{C(y,\P^{\Psi} \cup \{ \0 \} )}
|w-y|^r dw \exp[-\Delta(\{ {\bf 0}, y \}, {\cal
   P}^{\Psi})]  \right]-
(\E M^{\Psi}(\tau) )^2 \right) dy. $$

\vskip.5cm For any random point measure $\rho$, recall that
$\overline{\rho}$ denotes its centered version, that is
$\overline{\rho} := \rho - \E \rho$.

\begin{theo}\label{Gibbsthm} 
Assume that the density $h$ of $P$ is continuous on $Q_1$. We have
for each 
 $f \in \B(Q_1)$ \be \label{GibbsLLN} \lim_{\la \to \infty} \la^{r/d}
\langle f, \mu^{\Psi}_\la \rangle = \tau \E [M^{\Psi}(\tau)]
\int_{Q_1} h(x) f(x)  dx \ \ \text{in} \ L^2
\ee 
and \be \label{vargibbslim} \lim_{\la \to \infty} \la^{1 + 2r/d}
\Var [ \langle f, \mu^{\Psi}_\la \rangle ] = \tau V^{\Psi}(\tau)
\int_{Q_1} f^2(x) h^2(x) dx. \ee The finite-dimensional
distributions $\la^{-1/2 + r/d}(\langle f_1,
\overline{\mu}^{\Psi}_{\la} \rangle,...,\langle f_k,
\overline{\mu}^{\Psi}_{\la} \rangle), \ f_1,...,f_k \in \B(Q_1),$
of the random measures $(\la^{-1/2 + r/d}
\overline{\mu}^{\Psi}_{\la})$ converge as $\la \to \infty$  to
those of a mean zero Gaussian field with covariance kernel
 \be \label{gibbsCLT} (f_1,f_2) \mapsto \tau V^{\Psi}(\tau) \int_{Q_1} f_1(x) f_2(x) h^2(x)  dx, \  f_1, f_2  \in
\B(Q_1). \ee
\end{theo}

 \vskip.5cm


{\bf  Upper bounds for quantization coefficients.}  When $f \equiv
1$ the right hand side of (\ref{GibbsLLN}) gives
$$
\lim_{\la \to \infty} \la^{r/d} \langle 1, \mu^{\Psi}_\la \rangle
= \lim_{\la \to \infty} \la^{r/d} I(\tilde\P^{\Psi}_{(\la)})=
\tau \E [M^{\Psi}(\tau)].
$$
Since the  Bucklew and Wise limit (\ref{BW}) is necessarily no
larger than the right hand side of the above, this shows that in
addition to the bound (\ref{upbd}), that the $r$th quantization
coefficient $Q_{r,d}$ also satisfies the upper bound
$$
Q_{r,d} \leq (||h||_{d/(d+ r)})^{-1} \tau \E [M^{\Psi}(\tau)].
$$
Recall
 from our discussion above that when $\Psi \equiv 0$
(i.e. $\P^{\Psi}$ is Poisson) and when $f \equiv 1$, then the
right hand side of (\ref{GibbsLLN}) equals
$\tau^{-r/d} \omega_d^{-r/d} \Gamma(1 + r/d)$ and thus
$$
 Q_{r,d} \leq (||h||_{d/(d+ r)})^{-1} \tau^{-r/d} \omega_d^{-r/d}
 \Gamma(1 + r/d).
$$
We believe, although are not yet able to provide a full proof,
that whereas the distortion error (\ref{GibbsLLN}) is relatively
large for Poisson input, it can be made smaller if we restrict to
point sets which themselves enjoy some built-in repulsivity while
keeping the same mean point density. Indeed, given a fixed mean
number of test points it seems more economical to spread them
equidistantly over the domain of target distribution than to allow
for local overfulls of test points in some regions, which only
result in wasting test resources with the quantization quality
improvement considerably inferior to that which would be achieved
should we shift the extraneous points to regions of lower test
point concentration. In other words, the right hand side of
(\ref{GibbsLLN}) for repulsive Gibbs point processes should be
smaller than the corresponding distortion for the Poisson point
process with the same point density. It should be emphasized here
that in order to stay within the set-up of our asymptotic theory
we have to assume that the repulsivity is weak. On the other hand,
it is very likely that going to some extent beyond this
requirement may lead to even smaller quantization errors. These
seem to be natural and interesting questions, yet at present we
cannot handle them with our current techniques. 

\vskip.2cm

{\em Proof of Theorem \ref{Gibbsthm}.}  We claim that the
assertions of Theorem \ref{Gibbsthm} can be
  reduced to an application of Theorem \ref{BDDLIMIT} for functionals with bounded perturbations.
   We do it first assuming that the density
  $h$ is bounded away from $0.$ To this end, consider the following
  parametric family of geometric functionals:
  \begin{equation}\label{QUANTXIL}
   \hat\xi(x,\X;\la) := \int_{C(x,\X)} |y-x|^r
    \frac{P(d\la^{-1/d}y)}{\la h(\la^{-1/d}x)} =
    \int_{C(x,\X)} |y-x|^r \frac{h(\la^{-1/d}y)}{h(\la^{-1/d}x)} dy.
  \end{equation}
  Putting
  \begin{equation}\label{QUANTXI}
   \xi(x,\X) := \int_{C(x,\X)} |y-x|^r dy
  \end{equation}
  we obtain the bounded perturbed  representation (\ref{BDDIST1})
  for $\hat\xi(\cdot,\cdot;\la)$ with
  \begin{equation}\label{QUANTDELTA}
   \delta(x,\X;\la) := \int_{C(x,\X)} |y-x|^r \frac{h(\la^{-1/d}y)-
    h(\la^{-1/d}x)}{h(\la^{-1/d}x)} dy.
  \end{equation}
  It is easily seen that on \sr  \ input both $\xi$ and $\delta$ as given
  in (\ref{QUANTXI}) and (\ref{QUANTDELTA})  stabilize
  exponentially with common stabilization radius determined by
  the diameter of the Voronoi cell around the input point; see
   Subsection \ref{Euclid}(ii).
  We claim that this $\delta(\cdot,\cdot;\la)$  also satisfies the bounded
  moments condition (\ref{BDDIST2}). To this end, use the fact that
  $h$ is bounded away from $0$ and  write
  \begin{equation}\label{DeltaBound}
   |\delta(x,\X;\la)| \leq C
    \int_{C(x,\X)} |y-x|^{r} (h(\la^{-1/d}y)-
    h(\la^{-1/d}x))dy.
  \end{equation}
  Letting $\X := \P^{\Psi}$ and using the exponential
  decay of the Voronoi cell diameter on \sr \  input, see
  Subsection \ref{Euclid}, we conclude from (\ref{DeltaBound})
  that, for all $p > 0$ and $\la > 0$
  $$ \sup_x \E [\delta(x,\P^{\Psi};\la)]^p \leq
     \sup_x C \E \left[\int_{C(x,\P^{\Psi})} |y-x|^{r} (h(\la^{-1/d}y)-
     h(\la^{-1/d}x))dy \right]^p. $$
  By the translation invariance of $\P^{\Psi}$ and by the uniform continuity of
  the density $h$, this is bounded above by
  $$ C \E \left[\int_{C({\bf 0},\P^{\Psi})} |y|^r \omega_h(\la^{-1/d}|y|) dy\right]^p =: L(p,r,\la), $$
  where $\omega_h(\cdot)$ is the modulus of continuity of $h.$
  Since $\left[\int_{C({\bf 0},\P^{\Psi})} |y|^r \omega_h(\la^{-1/d}|y|) dy\right]^p$ is dominated
  by an integrable function of $\omega$ uniformly over $\la$, namely by a constant multiple of
  the $p(r+d)^{th}$ power of the Voronoi cell diameter,
  and since $\left[\int_{C({\bf 0},\P^{\Psi})} |y|^r \omega_h(\la^{-1/d}|y|) dy\right]^p$
  converges to zero as $\la \to \infty$ for almost all
  $\omega$,  we may use
  the bounded  convergence theorem to conclude that $L(p,r,\la) \to 0$ as $\la \to \infty$.
    This is clearly enough to get (\ref{BDDIST2}) and hence $\hat\xi(\cdot,\cdot;\la)$
  is an asymptotically negligible bounded perturbation of $\xi(\cdot,\cdot).$
  To proceed, we note that the quantization empirical measure
  $\mu^{\Psi}_{\la}$ satisfies for each $f \in \B(Q_1)$
  \begin{equation}\label{MPSI}
   \langle f, \mu^{\Psi}_{\la} \rangle = \la^{-1-r/d} \langle f h,
           \mu^{\hat\xi}_{\la} \rangle,
  \end{equation}
  where $\mu^{\hat\xi}_{\la}$ is the standard empirical measure
  (\ref{MUX}) for $\hat\xi$ as in (\ref{QUANTXI}), that is to
  say
  $$ \mu^{\hat\xi}_{\la} := \sum_{x \in \P^{\Psi} \cap Q_{\la}}
     \hat\xi(x,\P^{\Psi} \cap Q_{\la};\la) \delta_{x/\la^{1/d}}. $$
  On the other hand, it is easily verified that $\xi$ satisfies
  all assumptions of our limit Theorems \ref{WLLN}, \ref{VAR}
  and \ref{CLT}. Consequently, Theorem \ref{BDDLIMIT} can be
  applied for $\hat\xi,$ which yields
  Theorem \ref{Gibbsthm} via the formula (\ref{MPSI}) allowing
  us to translate results for $\mu^{\hat\xi}_{\la}$ to the
  corresponding results for $\mu^{\Psi}_{\la}.$ This  completes the proof of Theorem \ref{Gibbsthm}
  for $h$ bounded away from $0.$

 To proceed, assume now that $h$ fails to be bounded away from
  $0$ and, for $\varepsilon > 0$ put $h_{\varepsilon} :=
  \max(h,\varepsilon)$ and let $\mu^{\Psi}_{\la;\varepsilon}$
  be the version of $\mu^{\Psi}_{\la}$ with $h$
  replaced by $h_{\varepsilon}.$ Using the definition of
  $\mu^{\Psi}_{\la},$ and the exponential decay of the diameter
  of Voronoi cells in a Poisson-like  environment we easily conclude that
  \begin{equation}\label{EDIF}
     |\E[\langle f, \mu^{\Psi}_{\la} \rangle - \langle f, \mu^{\Psi}_{\la;\varepsilon} \rangle]|
     = O(\la^{-r/d} \varepsilon).
  \end{equation}
  Likewise, using the same we get
  \begin{equation}\label{VDIF}
   \Var[\langle f, \mu^{\Psi}_{\la} \rangle - \langle f, \mu^{\Psi}_{\la;\varepsilon} \rangle]
   = O(\la^{-1-2r/d} \varepsilon).
  \end{equation}
  Applying Theorem \ref{Gibbsthm} for $h_{\varepsilon},$ which is legitimate
  due to $h_{\varepsilon}$ being bounded away from $0,$ and then using
  (\ref{EDIF}) and (\ref{VDIF}) we readily get the required expectation
  and variance asymptotics for $\langle f, \mu^{\Psi}_{\la} \rangle$
  as well as the $L^2$ weak law of large numbers, which follows by
  the variance convergence. The remaining central limit theorem statement for
  $\langle f, \bar{\mu}^{\Psi}_{\la} \rangle$ follows directly by
  the Stein method as in Theorem \ref{CLT}, which is not affected
  by $h$ being not bounded away from $0.$ This completes the proof
  of Theorem \ref{Gibbsthm} for general $h.$ \qed

{\em Acknowledgements.}  J. E. Yukich thanks Ilya Molchanov for
helpful discussion related to quantization.

\vskip.5cm

Tomasz Schreiber, Faculty of Mathematics and Computer Science,
Nicholas Copernicus University, Toru\'n, Poland: \ {\texttt
tomeks@mat.uni.torun.pl }


J. E. Yukich, Department of Mathematics, Lehigh University,
Bethlehem PA 18015:
\\
{\texttt joseph.yukich@lehigh.edu}

\end{document}